\documentclass[12pt]{amsart}
\usepackage{graphics}

\newtheorem{theorem}{Theorem}
\newtheorem{lemma}[theorem]{Lemma}
\newtheorem{proposition}[theorem]{Proposition}

\newtheorem{corollary}[theorem]{Corollary}

\theoremstyle{definition}
\newtheorem{definition}[theorem]{Definition}
\newtheorem{notation}[theorem]{Notation}
\newtheorem{example}[theorem]{Example}

\theoremstyle{remark}
\newtheorem{remark}[theorem]{Remark}

\newcounter{figcount} 
\def\figno{\addtocounter{figcount}{1}\thefigcount}

\def\ab{\allowbreak}
\allowdisplaybreaks
\newcommand{\C}{{\mathbb{C}}}
\def\cl{{\rm cl}}
\newcommand{\D}{{\mathcal D}}
\def\ds{\displaystyle}
\def\E{\mathcal{E}}
\def\myemptyset{\mathchar"001F}
\renewcommand{\emptyset}{\myemptyset}
\newcommand{\FF}{\mathcal{F}}
\newcommand{\HH}{{\mathcal H}}
\newcommand{\la}{\langle}
\def\ncc{{\it NCC}}
\def\ncl{{\it NCL}}
\let\ol=\overline
\def\ot{\otimes}
\newcommand{\R}{ {\mathbb R} }
\newcommand{\ra}{\rangle}
\def\Tr{{\rm Tr}}
\def\wt{\widetilde}

\def\qed{{\unskip\nobreak\hfil\penalty50
 \hskip2em\hbox{}\nobreak\hfil %{\sc qed}
\vbox to 7.7pt{\parindent0pt%%%%%%%%
\hsize7.7pt\hrule\vrule height7.3pt%
\hfill\vrule height7.3pt\hrule}
\parfillskip=0pt \finalhyphendemerits=0 \par}}

\def\figbox#1{%
\hbox{\vrule\vbox to 150pt{\hsize150pt\hrule
\vfill \noindent
\hbox to 150pt{\hfill\vbox{\hsize125pt\raggedright\noindent
#1}\hfill}\vfill
\hrule}\vrule}}

\def\falsebox#1{%
\vbox to 50pt{\parindent0pt\hrule
\hbox to 50pt{\vrule height 25pt depth 25pt
\hfill #1 \hfill\vrule}\hrule}}

\newcommand{\myatop}[2]{\genfrac{}{}{0pt}{}{#1}{#2}}

\title[Orthogonal Polynomials]
{Orthogonal Polynomials and \\
Fluctuations of Random Matrices}

\author[T. Kusalik]{Timothy Kusalik $^{(\ddagger)}$}

\address{Queen's University, Department of Mathematics and
  Statistics, Jeffery Hall, Kingston, ON K7L 3N6, Canada}

\thanks{$^\ddagger$ Research supported by a USRA from the
Natural Sciences and Engineering Research Council of Canada}

\email{2tpk@qlink.queensu.ca}

\author[J. A. Mingo]{James A. Mingo $^{(*)}$}

\thanks{$^*$ Research supported by Discovery Grants and a Leadership
Support Initiative Award from the Natural Sciences and Engineering
Research Council of Canada}

\email{mingo@mast.queensu.ca}

\author[R. Speicher]{\hbox{Roland Speicher $^{(*)(\dagger)}$}}

\thanks{$^\dagger$ Research supported by a Premier's
  Research Excellence Award from the Province of Ontario}

\email{speicher@mast.queensu.ca}

\begin{document}

\begin{abstract}
In this paper we establish a connection between the
fluctuations of Wishart random matrices, shifted Chebyshev
polynomials, and planar diagrams whose linear span form a
basis for the irreducible representations of the annular
Temperly-Lieb algebra.
\end{abstract}

\maketitle

\section{Introduction}

Wishart matrices are a family of matrices studied in the
statistics literature since 1928.  Besides the Gaussian
random matrices they constitute the most important random
matrix ensemble.  They can be described as follows. Let
$G_{M,N}$ be a $M \times N$ matrix whose entries are
independent complex Gaussian random variables with mean 0
and complex variance $1/N$. Let $X_N = G_{M,N}^\ast
G_{M,N}$. $X_N$ is a complex Wishart matrix (of parameter
$c=M/N$).

The fundamental quantities of interest for random matrix
ensembles are the asymptotic eigenvalue distribution and the
fluctuations around this asymptotics.  Whereas the main
questions about eigenvalue distributions have been mostly
answered decades ago, investigations around fluctuations are
more recent and there is currently a lot of interest in this
topic, in particular, in connection with the question of
universality.  In the case of Wishart matrices, the large
$N$ limit of the eigenvalue distribution was found in 1967
by Marchenko and Pastur and is now named after them.  The
question of fluctuations was addressed for the first time by
Jonsson \cite{jon} in 1982 and, much more recently and much
more detailed, by Cabanal-Duvillard \cite{thierry} in 2001.

Before we say more about their results let us first describe
the general picture.  In the following we will, for better
legibility, systematically suppress the index $N$ (or $M$)
at our random matrices. It is to be understood that random
matrices are $N\times N$-matrices and asymptotic statements
refer to the limit $N\to\infty$.  For many random matrices
$Y$ the family of random variables $\{\Tr(Y^n)\}_n$ becomes
asymptotically, as the size $N$ of the matrices goes to
infinity, Gaussian. The two questions addressed above
consist then in understanding the limit of the expectation
and of the covariance of these Gaussian random
variables. For the latter one would in particular like to
diagonalize it. Whereas the expectation (i.e., the
eigenvalue distribution) depends on the considered ensemble,
the covariance (i.e., the fluctuations) seem to be much more
universal.  There are quite large classes of random matrices
which show the same fluctuations. The most important class
is the one which is represented by the Gaussian random
matrices. Its fluctuations are diagonalized by Chebyshev
polynomials, see Johansson \cite{kurt}.

In the case of the Wishart matrices $X$, the asymptotic
Gaussianity of the traces was shown by Jonsson; the explicit
form of the covariance, however, was revealed only recently
by Cabanal-Duvillard \cite{thierry}.  He found polynomials
$\{ \Gamma_n \}_n$, which were shown to be shifted Chebyshev
polynomials, such that the random variables $\{ \Tr(
\Gamma_n(X) \}_n$ are asymptotically Gaussian and
independent in the large $N$ limit; that is the polynomials
$\{ \Gamma_n \}_n$ diagonalize asymptotically the
covariance.  Cabanal-Duvillard's approach relies heavily on
stochastic calculus.  In this paper we want to give a
combinatorial proof of his results which rests on a
combinatorial interpretation for the polynomials
$\Gamma_n$. This combinatorial approach allows very
canonically an extension of Cabanal-Duvillard's results to a
family of independent Wishart matrices, yielding our main
result.

\medskip\noindent{\bf
Theorem}\ {\it
Let $\{\Gamma_n\}_n$ be the shifted Chebyshev polynomials of
the first kind as considered by Cabanal-Duvillard and let
$\{ \Pi_n \}_n$ be the orthogonal polynomials of the
Marchenko-Pastur distribution (which are shifted Chebyshev
polynomials of the second kind).  Let $X_1, \dots , X_p$ be
independent Wishart matrices and consider in addition to
$\Tr(\Gamma_n(X_i))$ also, for $k\geq 2$, the collection of
random variables $\Tr( \Pi_{m_1}( X_{i_1}) \cdots
\Pi_{m_k}(X_{i_k}) )$, where the Wishart matrices which appear
must be cyclically alternating, i.e.,
$i_1\not=i_2\not=i_3\not=\dots\not=i_k \not=i_1$.  These
latter traces depend only on the equivalence class of $(i_1,
\dots , i_k)$ and $( m_1, \dots m_k)$ under cyclic
permutation. Assuming that we have chosen one representative
from each equivalence class, the random variables
$$\{ \Tr( \Gamma_n( X_i ) ) \}
\cup
%\qquad\text{and}\qquad 
\{ \Tr(
\Pi_{m_1}( X_{i_1}) \cdots \Pi_{m_k}(X_{i_k}) ) \}$$ are
asymptotically independent and Gaussian.}

Let us remark that, in contrast to one-matrix models,
multi-matrix models are very poorly understood and the
problem of universality is mostly open for them at the
moment. Understanding the fluctuations of the simplest
representatives of multi-matrix models is essential for
progress in this direction.  Whereas our results about the
multi-matrix Wishart case are new, the corresponding results
for the multi-matrix Gaussian case were derived by
Cabanal-Duvillard in \cite{thierry}. By a small modification
of our approach we can also give a combinatorial
re-interpretation and proof in this case (see the remarks in
Section 11).

The main motivation for our investigations comes from our
belief that the theory of free probability provides the
right tools and concepts for attacking questions on
fluctuations of random matrices - in particular, for
multi-matrix models.  Even though ``freeness" will not appear
explicitly in this paper, our methods and results are very
much related to our investigations around ``second order
freeness" in \cite{ms}. The present paper can, in
particular, be seen as a complementary treatment of some of
the questions treated in \cite{ms}.

In the rest of this Introduction we want to give some idea
of what is involved in the proof of our theorem; in
particular, we would like to outline the relation between
special planar diagrams and the question of diagonalizing
the covariance of Wishart matrices.

Let us start with our re-interpretation of
Cabanal-Duvillard's results for the case of one Wishart
matrix.  Our starting point is the paper of Mingo-Nica
\cite{mn}, where a genus expansion in terms of permutations
was provided for the cumulants of the random variables
$\Tr(X^n)$. Since cumulants of different orders have
different leading contributions in $N$, this has as a direct
consequence the asymptotic Gaussianity of these traces.  The
main problem left is to understand and diagonalize the
covariance. Also in \cite{mn}, it was shown that the
covariance of the random variables $\{ \Tr( X^n )\}_n$ has
asymptotically a very nice combinatorial interpretation,
namely it is given by counting a class of planar diagrams
which were called non-crossing annular permutations.  More
precisely, if we denote by $\kappa_2(A,B)$ the covariance of
two random variables $A$ and $B$ and if $c$ is the
asymptotic ratio of $M$ and $N$ for our Wishart matrices, we
have 
$$\lim_{N\to\infty}
\kappa_2(\Tr(X^m),\Tr(X^n))=\sum_{\pi\in S_{NC}(m,n)}
c^{\#\pi},$$ 
where $S_{NC}(m,n)$ denotes the set of non-crossing
$(m,n)$-annular permutations, i.e., permutations on $m+n$
points which connect $m$ points on one circle with $n$
points on another circle in a planar or non-crossing way. In
the above formula we are summing over all non-crossing
$(m,n)$-annular permutations and each block of such a
permutation contributes a multiplicative factor $c$. (For
$c=1$, which corresponds to Wishart matrices with $M=N$, the
above formula just counts the number of elements in
$S_{NC}(m,n)$.)

In \cite{mn} the cycles of the permutation were shown as
blocks in the annulus. See the figure below.

\begin{center}
$\vcenter{\hsize175pt\includegraphics{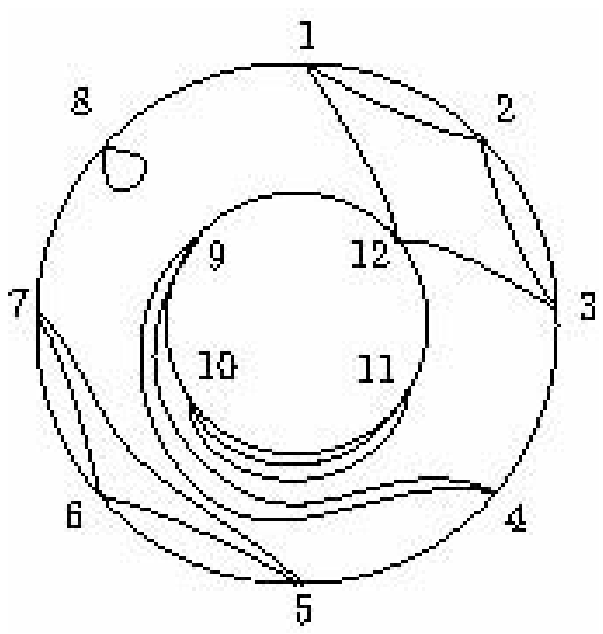}}$\qquad
$\vcenter{\hsize 150pt\raggedright\small
Figure \figno. A non-crossing annular permutation.}$
\end{center}

In our context it seems more appropriate to redraw the
diagram with the circles side by side as shown. 

\begin{center}
$\vcenter{\hsize255pt\includegraphics{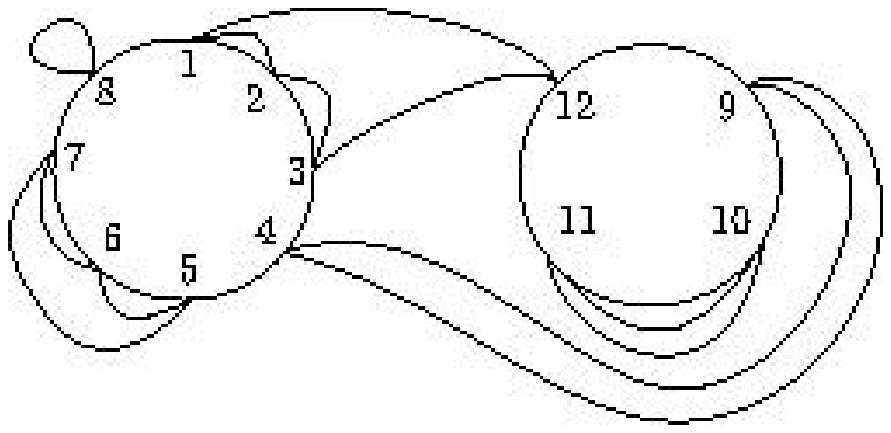}}$

\
\vbox{{\small Figure \figno. The permutation of figure 1
unfolded.}}
\end{center}

The limit of $\kappa_2(\Tr(X^m),\Tr(X^n))$ is, of course,
not diagonal in $m$ and $n$ because points on each circle
are grouped into blocks (some of which do not even connect
to the other circle), and this grouping on both sides has no
correlation; so there is no constraint that $m$ has to be
equal to $n$.  However, a quantity which clearly must be the
same for both circles is the number of through-blocks, i.e.,
blocks which connect both circles.  Thus in order to
diagonalize the covariance we should go over from the number
of points on a circle to the number of through-blocks
leaving this circle.

A nice way to achieve this is to cut our diagrams in two
parts -- one part for each circle as shown.

\begin{center}
\includegraphics{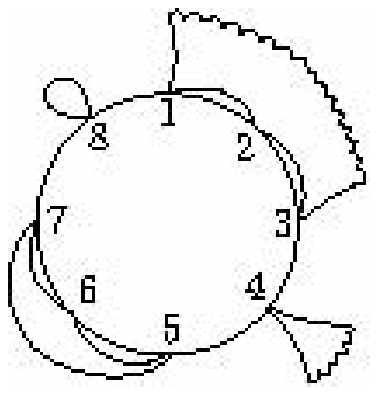}\qquad\qquad
\includegraphics{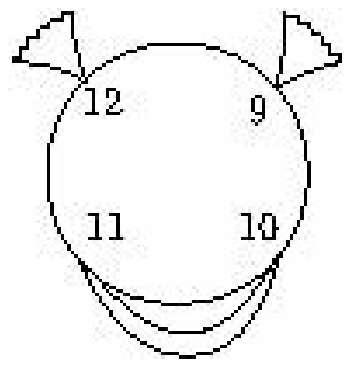}
\vbox{\small\textbf{Figure \figno.} The permutation of figure 2
cut in two.} \end{center}

These diagrams which we call `non-crossing
half-permutations' are the principal objects of study of
this paper.

In this pictorial description $\Tr(X^m)$ corresponds to the
sum over non-crossing half-permutations on one circle with
$m$ points and $\Tr(X^n)$ corresponds to a sum over
non-crossing half-permutations on another circle with $n$
points. The limit of $\kappa_2(\Tr(X^m),\Tr(X^n))$
corresponds to pairing the half-permutations for $\Tr(X^m)$
with the half-permutations for $\Tr(X^n)$. The pairing
between two half-permutations is given by glueing them
together in all possible planar ways.  We will call what is
left in a half-permutation of a through-block after cutting
an `open' block -- as opposed to `closed' blocks which live
totally on one circle and are thus not affected by the
cutting. A pairing of two half-permutations consists of
glueing together their open blocks in all possible planar
ways. This clearly means that both half-permutations must
have the same number of open blocks for a non-trivial
pairing, and thus our covariance should become diagonal if
we go over from the number of points on a circle to the
number of open blocks.  From this point of view the
polynomials $\Gamma_k(x)$ found by Cabanal-Duvillard should
describe $k$ open blocks.  If we write $x^n$ as a linear
combination of the $\Gamma_k$'s as in
$$x^n = \sum_{k=0}^n q_{n,k} \Gamma_k(x)$$
then the above correspondence
\begin{align*}
\Tr(X^n)\quad&\hat=\quad\text{$n$ points}\\
\Tr(\Gamma_k(X))\quad&\hat=\quad\text{$k$ open blocks}
\end{align*}
suggests that the coefficients $q_{n,k}$ are given by
summing over all half-permutations with $k$ open blocks,
each such permutation contributing a factor $c$ for each
closed block.  (The dependence on $c$ reflects the fact that
in our original formula every block contributed with a
factor $c$ -- now every closed block gives a factor $c$
right away, whereas an open block has to be paired with
another open block to produce a factor $c$. In the case
$c=1$, $q_{n,k}$ just counts the number of half-permutations
on $n$ points which have $k$ open blocks.)  We will show
that the coefficients $q_{n,k}$ appearing in the relations
for the shifted Chebyshev polynomials $\Gamma_k$ have indeed
this combinatorial meaning. This will be achieved, in
Section 8, by showing that both quantities satisfy the same
recurrence relations. Let us also point out that the case
$k=0$ is special, because constant terms in the polynomials
play no role for the covariance, but have to be fixed by
other normalizations. On the combinatorial level this is
reflected by the fact that we only look on non-crossing
annular permutations which connect the two circles, thus we
always have at least $k\geq 1$ through-blocks (or open
blocks, after cutting). Since however, our recursions rely
also on $k=0$, we have to make some separate considerations
for $k=0$ at various places. In particular, we want to point
out that the `right' definition in our setting for a
non-crossing half-permutation on $n$ points with zero open
block is \emph{not} just a non-crossing permutation on $n$
points. See Section \ref{recursion} for more details on this.

In order to illustrate the above statements let us here present
the pictorial meaning of the equation
$$x^2=\Gamma_2+(2+2c)\Gamma_1+q_{2,0}\Gamma_0$$
for the shifted Chebyshev polynomials $\Gamma_k$.

$$ x^2 =
\mathop{\vcenter{\hbox{\includegraphics{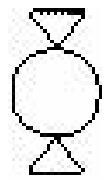}}}}_{\hbox{\vrule
width 0pt height 15pt$\Gamma_2$}}
\quad + \quad
\mathop{\vcenter{\hbox{\includegraphics{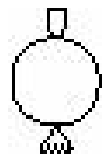}}}}_{\hbox{\vrule
width 0pt height 15pt$c \Gamma_1$}}
\quad + \quad
\mathop{\vcenter{\hbox{\includegraphics{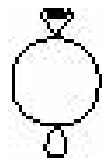}}}}_{\hbox{\vrule
width 0pt height 15pt$c \Gamma_1$}}
\quad + \quad
\mathop{\vcenter{\hbox{\includegraphics{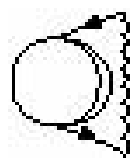}}}}_{\hbox{\vrule
width 0pt height 15pt$\Gamma_1$}}
\quad + \quad
\mathop{\vcenter{\hbox{\includegraphics{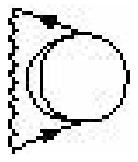}}}}_{\hbox{\vrule
width 0pt height 15pt$\Gamma_1$}}
$$ 
{\small\begin{center}\leavevmode
\vbox{\hsize 285pt\raggedright\small
\textbf{Figure \figno.} The
five non-crossing half-permutations on $[2]$. The
diagrams corresponding to $\Gamma_0$ have been
omitted.}\end{center}}

Our combinatorial interpretation of the diagonalizing
polynomials for one Wishart matrix opens the way to a
similar treatment for a tuple of independent Wishart
matrices. Covariances of traces of products of such matrices
are asymptotically described essentially in the same way as
before in terms of non-crossing annular permutations,
however, in addition we have to colour the points on the two
circles (one colour for each Wishart matrix) and require
that the contributing permutations have to connect only
points of the same colour (i.e., each cycle of the
permutation must be mono-chromatic). Again the
diagonalization of the covariance is achieved by going over
from the number of points on a circle to the number of open
blocks. Thus, on first view, one might expect that traces of
alternating products in the $\Gamma_k$ give rise to a
diagonal covariance. However, this is \emph{not} the case.
One has to realize that through-blocks of one colour break
the symmetry of the circle for through-blocks of another
colour, thus in a sense for each group of through-blocks of
the same colour the circle is cut open to a line and instead
of circular half-permutations we have to consider linear
half-permutations. Thus instead of $\Gamma_k$ we have to
look for polynomials $\Pi_k$ which satisfy
$$x^n=\sum_{k=0}^n p_{n,k} \Pi_k(x),$$
where $p_{n,k}$ is now the sum over all linear
half-permutations with $k$ open blocks, each closed block
weighted by a factor $c$.  It turns out that these
polynomials are the orthogonal polynomials for the
Marchenko-Pastur distribution (which are shifted Chebyshev
polynomials of the second kind).  Again we prove this, in
Section 7, by showing the equality of the corresponding
recurrence relations. The proof that the covariance is
diagonalized by traces in alternating products in these
polynomials $\Pi_k$ is given in Section 9.

As an illustrative example for these statements consider 
$$
x^2 y = \big\{\Pi_2(x) + (1 + 2 c) \Pi_1(x) + (c + c^2) \Pi_0(x)
\big\}\, \big\{\Pi_1(y) + c \Pi_0(y)\big\}
$$
 
$$ x^2 y = \kern-10pt
\mathop{\vcenter{\hbox{\includegraphics{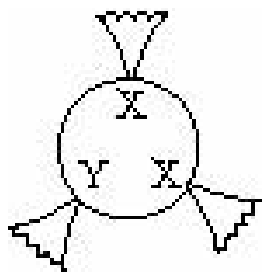}}}}_{\hbox{\vrule
width 0pt height 15pt$\Pi_2(X) \Pi_1(Y)$}}
 + 
\mathop{\vcenter{\hbox{\includegraphics{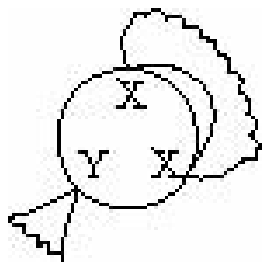}}}}_{\hbox{\vrule
width 0pt height 15pt$\Pi_1(X) \Pi_1(Y)$}}
\ + 
\mathop{\vcenter{\hbox{\includegraphics{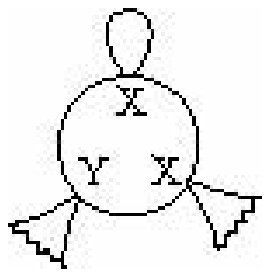}}}}_{\hbox{\vrule
width 0pt height 15pt$c \Pi_1(X) \Pi_1(Y)$}}
 + 
\mathop{\vcenter{\hbox{\includegraphics{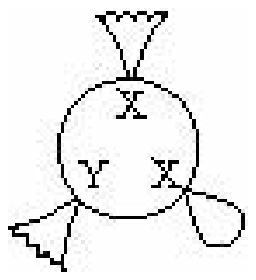}}}}_{\hbox{\vrule
width 0pt height 15pt$c \Pi_1(X) \Pi_1(Y)$}}
$$ 
{\small\begin{center}\leavevmode
\vbox{\hsize 300pt\raggedright\small
\textbf{Figure \figno.} The
four non-crossing circular half-permutations on $\{X, X, Y\}$.
The diagrams corresponding to $\Pi_0$ have been
omitted.}\end{center}}

Note that the above argument of symmetry breaking is strictly
valid only if each group of the same colour has at least one
through-block. However, a priori one also has to consider
diagrams containing mono-chromatic groups without any
through-block. It turns out that, by the centeredness of the
$\Pi_k$, the net contribution of such diagrams cancels out. To
make this argument rigourous constitutes an essential part of
the proof of Theorem \ref{mt} (see in particular Lemma
\ref{lemma18})

Finally, we would like to point out that our circular
half-permutations are, after a small redrawing, the diagrams
used by V. F. R.~Jones \cite{jones}, \S 5 to create a basis
for the irreducible representations of the annular
Temperly-Lieb algebras. In Jones's convention the left
picture in Fig. 3  above would be first inverted in the
centre of the circle and then the blocks would be expanded
into `fat graphs'.

\begin{center}
\includegraphics{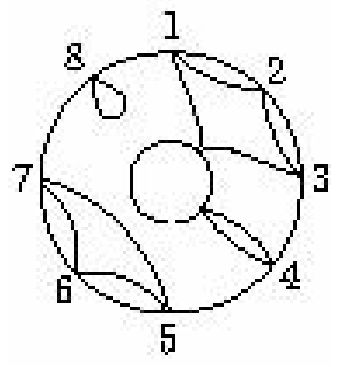}\qquad\qquad
\includegraphics{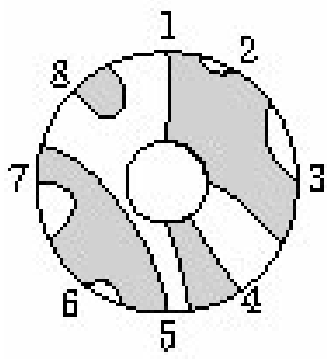}
{\small\leavevmode
\vbox{\hsize 285pt\raggedright\small
\textbf{Figure \figno.} 
The left diagram of figure 3 inverted and then drawn as a fat
permutation. }} \end{center}

In this paper we will not have to say more about this
relation with the Temperly-Lieb algebras, but we are taking
this as a serious hint that there exists a deeper relation
between free probability, random matrices, and subfactors.
We hope to explore this relation further in forthcoming
investigations.

However, as a nice application of our developed machinery
to operator algebraic questions we show, in Section 10, a
connection between Wick products and half-permutations that
gives a combinatorial formula for the product of two Wick
products which has a very simple diagrammatic
interpretation.

\begin{eqnarray*}
\raise-4pt\hbox{\includegraphics{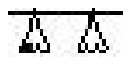}\ } \ast 
\raise-4pt\hbox{\ \includegraphics{fig20.eps}} & =&
\raise-4pt\hbox{\includegraphics{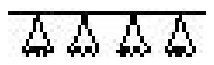}\ } +
\raise-4pt\hbox{\ \includegraphics{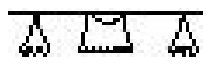}\ } +
\raise-4pt\hbox{\ \includegraphics{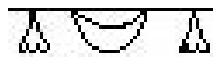}} \\ [10pt]
&& \mbox{} +
\raise-4pt\hbox{\ \includegraphics{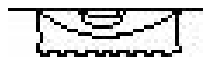}\ } +
\raise-4pt\hbox{\ \includegraphics{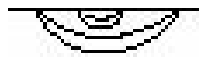}}
\end{eqnarray*}
\begin{center}\leavevmode{\small
\vbox{\hsize 250pt\raggedright
\textbf{Figure \figno.} The convolution of two non-crossing
linear half permutations.}}
\end{center}

\section{Notation}

Let $(\Omega, \Sigma, {\sf P})$ be a probability space. We
will only consider random variables with moments of all
orders. Let $G_{M, N}: \Omega \rightarrow M_{M, N}({\mathbb
C})$ be a random matrix with entries $(g_{i,j})$ such that
each $g_{i,j}$ is a complex Gaussian random variable with
mean 0 and complex variance $\E(|g_{i,j}|^2) = 1/N$. Let
$X_N = G_{M, N}^\ast G_{M, N}$, we shall call $X_N$ a
complex Wishart distributed random matrix or more succinctly
a {\em Wishart matrix}. See for example Haagerup and
Thorbj{\o}rnsen \cite{ht1}, \S 5, \cite{ht2}, \S 6 or
Hiai and Petz \cite{hp}, Ch. 4.

For any random variables {\sf X} and {\sf Y}, let the covariance
be given by  $\kappa_2({\sf X}, {\sf Y})\ab = \E( {\sf X} {\sf
Y}) - \E({\sf X}) \E({\sf Y})$. For an $N \times N$ matrix $X$
we let $\Tr(X)$ be the un-normalized trace.

In this paper we are interested in the large $N$ limit of
the eigenvalue distribution of $X_N$ and their fluctuations. In order
for this limit to exist we must have that $M/N$ converges as $M$ and
$N$ tend to infinity.  More precisely we shall suppose that
we have sequences $\{M_k \}_k$ and $\{ N_k \}_k$ such that
$\lim_k M_k = \lim_k N_k = \infty$ and $\lim_k M_k/N_k = c$
with $ 0 < c < \infty$. In order to have a second order
distribution we require in addition the existence of a
second order limit $c' = \lim_k M_k - c N_k$. We shall
assume that the numbers $c$ and $c'$ are fixed throughout
the paper.

This limiting distribution is called the Marchenko-Pastur
distribution, $$\mu_c = (1 - c) \delta_0 + \frac{ \sqrt{ (b
- x)(x -a)} }{ 2 \pi x} dx$$ where $a = (\sqrt{c} - 1)^2$ ,
$b = (\sqrt{c} +1)^2$ , and $\delta_0$ is an atom at 0 which
disappears if $c \geq 1$. The absolutely continuous part $
\sqrt{ (b - x)(x -a)}/( 2 \pi x)$ has support equal to the
interval $[a, b]$.

\section{Shifted Chebyshev Polynomials of the First
Kind}\label{firstkind}

We shall make use of two families $\{ \Gamma_k(x) \}_k$ and
$\{ \Pi_k(x) \}_k$ of {\em shifted} Chebyshev
polynomials \textit{c.f.} \cite{thierry}. Let $\{ T_k(x) \}_k$ be
the Chebyshev polynomials of the first kind i.e. $T_k( \cos
\theta ) = \cos ( k \theta)$. Rescale these polynomials by
letting $C_0(x)= 1$ and for $n \geq 1$, $C_k(x) = 2 T_k(x/2)$.
Then $\{ C_k(x) \}_k$ are monic and orthogonal for the dilated
arc-sine law $\ds\frac{dx }{\pi \sqrt{4 - x^2}}$ on $[-2,
2]$.

We have
\vtop{\hsize200pt
\begin{align*}
C_0(x) &= 1 \\
C_1(x) &= x \\
C_2(x) &= x^2 - 2 \\
C_3(x) &= x^3 - 3 x \mbox{ and }\\
x C_n(x) &= C_{n+1}(x) + C_{n-1}(x) \mbox{ for } n > 1 
\end{align*}}

Now let us shift the sequence to produce our first sequence
$\{ \wt \Gamma_n(x) \}_n$. Let
$$ u = \frac{ x - ( 1 + c)}{\sqrt{c}}\quad \mbox{ then }\quad
\frac{du}{\sqrt{4 - u^2}} = \frac{dx}{\sqrt{(b - x)( x - a)}}$$
where $a = (\sqrt{c} - 1)^2$ and $b = (\sqrt{c} + 1)^2$. 

Now let 
$$\wt\Gamma_n(x) = \sqrt{c^n}\ C_n\left(\frac{x - ( 1 + c)}{\sqrt{c}}\right)
$$

Then 
\begin{align*}
\wt\Gamma_0(x) &= 1 \\
\wt\Gamma_1(x) &= x - (1 + c) \\
\wt\Gamma_2(x) &= x^2 - 2(1 + c) x + 1 + c^2 \\
\wt\Gamma_3(x) &= x^3 - 3 (1 + c) x^2 + 3(1 + c + c^2) x  -(1 + c^3)\\
x \wt\Gamma_n(x) &= \wt\Gamma_{n+1}(x) + (1 + c) \wt\Gamma_n(x) + c
\wt\Gamma_{n-1}(x) \mathrm{\ for\ } n = 0, 2, 3, 4, \dots
\end{align*} 
and $\{ \wt\Gamma_n(x) \}_n$ is a sequence of monic pairwise orthogonal
polynomials for the shifted arc-sine law 
$\ds\frac{dx}{\pi \sqrt{(b - x)( x - a)}}$ on $[a, b]$. 

We shall write $\wt\Gamma_n(x) = \sum_{k=0}^n g_{n,k}' x^k$. The recurrence
formula for the $\wt\Gamma_n$'s means that for $n = 0, 2, 3, 4, \dots$
we have 
\begin{equation}\label{geq1}
g_{n,k-1}' = g_{n+1,k}' + (1 + c) g_{n,k}' + c g_{n-1,k}'
\end{equation} and for $n = 1$
\begin{equation}\label{geq2}
g_{1,k-1}' = g_{2,k}' + (1 + c) g_{1,k}' + 2 c g_{0,k}
\end{equation}
where we have used the convention that $g_{n,k}' = 0$ if either $n < 0$ or $k <
0$, or $k > n$. Let $\wt\Gamma$ denote the upper-triangular matrix

\begin{equation*}
\wt\Gamma = 
\begin{pmatrix}
g_{0,0}' &          &          &          & \\
g_{1,0}' & g_{1,1}' &          &          & \\
g_{2,0}' & g_{2,1}' & g_{2,2}' &          & \\
g_{3,0}' & g_{3,1}' & g_{3,2}' & g_{3,3}' & \\
         & \ddots   & \ddots   &  \ddots  & \ddots \\
\end{pmatrix}
\end{equation*}

Let $S$ be the unilateral shift and $C = c S + cI + SS^\ast + S^\ast$
\begin{equation}\label{csdef} S = \begin{pmatrix}
0 & 0 &   & \\
1 & 0 & 0 & \\
  & 1 & 0 &  \ddots\\ 
  &   & \ddots & \ddots\\
\end{pmatrix}, \ 
C =
\begin{pmatrix}
c & 1     &        &         \\
c & 1 + c &   1    &         \\
  &   c   & 1 + c  & \ddots  \\
  &       & \ddots & \ddots   \\
\end{pmatrix}
 \end{equation}

Equations (\ref{geq1}) and (\ref{geq2}) imply that 
$$ \wt\Gamma  S^\ast  = C\wt\Gamma  
+ (I + c S)(I - S S^\ast) \wt\Gamma$$ 
After multiplying on the left and the right by
$\wt\Gamma^{-1}$ we have 
\begin{equation}\label{gcommutator}
\wt\Gamma^{-1} C = S^\ast\wt\Gamma^{-1}  - \wt\Gamma^{-1} (I + c
S)(I - S S^\ast)
\end{equation}

If we write 
\begin{equation*} \wt\Gamma^{-1} = 
\begin{pmatrix}
g_{0,0} &         &         &         &      \\
g_{1,0} & g_{1,1} &         &         &       \\
g_{2,0} & g_{2,1} & g_{2,2} &         &        \\
g_{3,0} & g_{3,1} & g_{3,2} & g_{3,3} &        \\
        & \ddots  & \ddots  &  \ddots & \ddots \\
\end{pmatrix}
\end{equation*}
then equation (\ref{gcommutator}) implies that for $k > 0$
\begin{equation}\label{geq3}
g_{n+1, k} = g_{n,k-1} + (1 + c) g_{n,k} + c g_{n,k+1}
\end{equation}
and for $k = 0$
\begin{equation}\label{geq4}
g_{n+1,0} = c g_{n, 1} + (1 + c) g_{n,0} +  c g_{n,1}
\end{equation}

By direct calculation we obtain that the first five rows of $\wt\Gamma^{-1}$ are
{\footnotesize\begin{equation*}\begin{pmatrix}
1                & 0                 & 0        & 0 & 0 \\  
1 + c                & 1                 & 0        & 0 & 0 \\ 
1 + 4c + c^2          & 2 + 2\,c          & 1        & 0 & 0 \\ 
1 + 9\,c + 9\,c^2 + c^3 & 3 + 9\,c + 3\,c^2 & 3 + 3\,c & 1 & 0 \\ 
1 + 16c + 36\,c^2 + 16\,c^3 + c^4 & 4 + 24\,c + 24\,c^2 + 4\,c^3 & 
6 + 16\,c + 6\,c^2 & 4 + 4\,c & 1\\
\end{pmatrix}\end{equation*}}

Let $d_0 = -1$, $d_1 = 1$ and $d_n = (-1)^n (c - 1)$ for $n
> 1$. Let 
$$\Gamma_n(x) = \wt\Gamma_n(x) + d_n$$ 
then
$\int_{\R} \Gamma_n(x) \, d\mu_c(x) = 0$ for $n = 0, 1,
2$. By equation \ref{gamma_pi} above, 
$$\int_{\R} \wt\Gamma_n(x) \, d\mu_c(x) = (-1)^{n}( 1 -c)$$
for $n > 1$ and thus $\int_{\R} \Gamma_n(x) \, d\mu_c(x) = 0$ for all
$n$.

Let us write $\Gamma_n(x) = \sum_{k=0}^n q_{n,k}' x^k$ and $x^n =
\sum_{k=0}^n q_{n,k} \Gamma_k(x)$. Let $\Gamma$ be the matrix
\begin{equation*}
\Gamma = \begin{pmatrix}
q_{0,0}' &          &          &               \\
q_{1,0}' & q_{1,1}' &          &                \\
q_{2,0}' & q_{2,1}' & q_{2,2}' &                 \\
q_{3,0}' & q_{3,1}' & q_{3,2}' & q_{3,3}' &       \\
         & \ddots   & \ddots   & \ddots   & \ddots \\
\end{pmatrix}
\end{equation*}and thus \begin{equation*}
\Gamma^{-1} = \begin{pmatrix}
q_{0,0} &         &         &               \\
q_{1,0} & q_{1,1} &         &                \\
q_{2,0} & q_{2,1} & q_{2,2} &                 \\
q_{3,0} & q_{3,1} & q_{3,2} & q_{3,3} &       \\
        & \ddots  & \ddots  & \ddots  & \ddots \\
\end{pmatrix}
\end{equation*}

Let $D$ be the matrix
\[ D =
\left(\begin{tabular}{r|cc}
$d_0$    & 0      & $\cdots$  \\
$d_1$    & 0      & $\cdots$  \\
$d_2$    & 0      & $\cdots$  \\
$d_3$    & 0      & $\cdots$  \\
$\vdots\ $ & \vdots &   \\
\end{tabular}\right)
\]

Then $\Gamma = \wt\Gamma + D$ and so $\Gamma^{-1} = \wt\Gamma^{-1} -
\wt\Gamma^{-1}D$, since $D \Gamma = D$ which holds because $q_{0,0}'
= 1$. Thus $q_{n,k} = g_{n,k}$ for $k > 0$. 

 The first five rows of $\Gamma^{-1}$ are
{\footnotesize\begin{equation*}\begin{pmatrix}
1                & 0                 & 0        & 0 & 0 \\  
c                & 1                 & 0        & 0 & 0 \\ 
c + c^2          & 2 + 2\,c          & 1        & 0 & 0 \\ 
c + 3\,c^2 + c^3 & 3 + 9\,c + 3\,c^2 & 3 + 3\,c & 1 & 0 \\ 
c + 6\,c^2 + 6\,c^3 + c^4 & 4 + 24\,c + 24\,c^2 + 4\,c^3 & 
6 + 16\,c + 6\,c^2 & 4 + 4\,c & 1\\
\end{pmatrix}\end{equation*}}
We shall see in section \ref{recursion} that these equations have a
very simple combinatorial interpretation.

Recall that the Cauchy transform of the arc-sine law is
$$
G(z) = \ds\int_{-2}^2 \frac{1}{z - t} \frac{dt}{\pi \sqrt{4
- t^2}} = \frac{1}{\sqrt{z^2 - 4}}
$$  
Using the
transformation $u = (z - (1 + c))/\sqrt{c}$ we see
that the Cauchy transform of the shifted arc-sine law is
$$
F(z) = \ds\int_a^b \frac{1}{z - t} \frac{dt}{\pi\sqrt{(b
-t)(t-a)}} = G(u)/\sqrt{c}
$$
Thus the moment generating
function of the shifted arc-sine law is 
$$
\sum_{n \geq 0}
g_{n,0} z^n = z^{-1} F(z^{-1}) = \ds\frac{1}{\sqrt{(1 -
az)(1 -bz)}}
$$

Now for each $k > 0$ let 
$$G_k(z) = \sum_{n \geq k} g_{n,k}
z^n$$ 
Then $G_0(z) = \frac{1}{\sqrt{(1 - az)(1 -
bz)}}$. From equation (\ref{geq3}) we see that for $k \geq
1$
\begin{equation}\label{Grec1}
c G_{k+1}(z) = (z^{-1} - (1 +  c)) G_k(z) -
G_{k-1}(z)
\end{equation}
and 
\begin{equation}\label{Grec2}
2 c G_1(z) = (z^{-1} - (1 + c)) G_0(z) - z^{-1}
\end{equation}
 
Let \[P_0(z) = \frac{1 - (c -1) z -\sqrt{(1 -az)(1 - bz)}}{2
z} \] be the moment generating function of the
Marchenko-Pastur distribution.  Then we can write equation
(\ref{Grec2}) as 
\begin{equation*} 
G_1(z) = \bigg( \frac{P_0(z) - 1}{c}\bigg) G_0(z) 
\end{equation*}
Moreover $P_0(z)$ satisfies the functional equation $(z^{-1} - (1 +
c)) (P_0(z) -1) = (P_0(z) -1)^2 + c$. Now we can apply
equation (\ref{Grec1}) to conclude by induction that
\begin{equation}\label{power1} G_n(z) = \bigg(\frac{P_0(z)
-1}{c}\bigg)^n G_0(z) \end{equation}

There is an interesting diagrammatic
interpretation of equation (\ref{power1}) in
section \ref{circular}, c.f. Remark \ref{circulargenfunction}.

\section{Shifted Chebyshev Polynomials of the Second
Kind}\label{pi-poly}

We next recall the construction of the orthogonal
polynomials for the Marchenko-Pastur distribution
$\mu_c$. Let $\{ U_n(x) \}_n$ be the Chebyshev polynomials
of the second kind, i.e. $U_n( \cos \theta) = \sin( (n+1)
\theta )/ \sin \theta$. Let $S_n(x) = U_n(x/2)$. Then

\vtop{\hsize200pt
\begin{align*}
S_0(x) &= 1 &  S_3(x) &= x^3 - 2 x \\
S_1(x) &= x & x S_n(x) &= S_{n+1}(x) + S_{n-1}(x) \mbox{ for } n > 1\\
S_2(x) &= x^2 - 1 \\
%S_3(x) &= x^3 - 2 x \mbox{ and }\\
% 
\end{align*}}
The $S_n$'s are monic and orthogonal for the semicircle law $\sqrt{4 - x^2}/(
2 \pi)$.  We let 
\begin{equation}\label{pi-def}
\Pi_n(x) = \sqrt{c^n}\ S_n\left( \frac{ x - (1 +
c)}{\sqrt{c} }\right) + \sqrt{c^{n-1} }\ S_{n-1}\left( \frac{ x - (1 +
c)}{\sqrt{c} }\right)
\end{equation}
Then 
\begin{align*}
\Pi_0(x) &= 1    &  \Pi_2(x) &= x^2 - (1 + 2 c) x + c^2\\
\Pi_1(x) &= x -c &  \Pi_3(x) &= x^3 - (2 + 3 c) x^2 + (1 + 2 c + 3c^2) x -
c^3\\
\end{align*} 
and for $n > 1$
$$
x \Pi_n(x) = \Pi_{n+1}(x) + (1 + c) \Pi_n(x) + c \Pi_{n-1}(x)
$$
The $\Pi_n$'s are shifted Chebyshev polynomial of the second kind.
Indeed by letting 
$$ u = \frac{ x - ( 1 + c)}{\sqrt{c}}\quad 
\mbox{ then }\quad 
\frac{c\,\sqrt{4 - u^2}\, du}{2 \pi} = \frac{\sqrt{(b - x)( x -
a)}\, dx}{2 \pi}$$
where $a = (\sqrt{c} - 1)^2$ and $b = (\sqrt{c} + 1)^2$, we
obtain, from equation (\ref{pi-def}), that $\int_{\R} x
\Pi_n(x) \, d\mu_c(x) = 0$ for $n > 1$. In addition we can
see by direct calculation that $\int_{\R} \Pi_1(x) \,
d\mu_c(x) = \int_{\R} \Pi_2(x) \, d\mu_c(x) = 0$. Thus by
the recurrence relation we have that $\int_{\R} \Pi_n(x) \,
d\mu_c(x) = 0$ for all $n > 0$.  Finally by induction on $m
+ n$ we have again from the recurrence relation for $n > m$ that
$\int_{\R} x^m \Pi_n(x)\, d \mu_c(x) = \int_{\R} x \Pi_{m+n
-1}(x) \, d \mu_c(x) = 0$. This shows that the $\Pi_n$'s are
pairwise orthogonal. From the recurrence relation we have that
$$
\int_{\R} \Pi_n(x)^2 \, d \mu_c(x) = \int_{\R} x \Pi_n(x)
\Pi_{n-1}(x)\, d \mu_c(x)
$$ 
and from equation (\ref{pi-def}) we obtain 
$$
\int_{\R} x \Pi_n(x) \Pi_{n-1}(x) \, d \mu_c(x) = c^n \int_{\R}
S_{n-1}(u)^2\, \frac{ \sqrt{ 4 - u^2}}{2 \pi} \, du = c^n
$$ 
Thus  $\int_{\R} \Pi_n(x)^2 \, d \mu_c(x) = c^n$.

Let us write $\Pi_n(x) = \sum_{k=0}^n p_{n,k}' x^k$ and $x^n
= \sum_{k=0}^n p_{n,k} \Pi_k(x)$. Let $\Pi$ be the matrix
\begin{equation*}
\Pi = \begin{pmatrix}
p_{0,0}' &          &          &               \\
p_{1,0}' & p_{1,1}' &          &                \\
p_{2,0}' & p_{2,1}' & p_{2,2}' &                 \\
p_{3,0}' & p_{3,1}' & p_{3,2}' & p_{3,3}' &       \\
         & \ddots   & \ddots   & \ddots   & \ddots \\
\end{pmatrix}
\end{equation*}and thus \begin{equation*}
\Pi^{-1} = \begin{pmatrix}
p_{0,0} &         &         &               \\
p_{1,0} & p_{1,1} &         &                \\
p_{2,0} & p_{2,1} & q_{2,2} &                 \\
p_{3,0} & p_{3,1} & p_{3,2} & q_{3,3} &       \\
        & \ddots  & \ddots  & \ddots  & \ddots \\
\end{pmatrix}
\end{equation*}

Let $C$ and $S$ be as in equation (\ref{csdef}) above. Then
we have from the recurrence relation for the $\Pi$'s that $C
\Pi = \Pi S^\ast$ and thus $\Pi^{-1}C = S^\ast
\Pi^{-1}$. From this we obtain the recurrence relation for
the $p_{n,k}$'s:
\begin{equation}\label{Prec1}
p_{n+1, k} = p_{n,k-1} + (1 + c) p_{n,k} + c p_{n, k+1} \mbox{,\ when\ } k > 0;
\end{equation}
\begin{equation}\label{Prec2}
\mbox{ and\ } p_{n+1,0} = c p_{n,0} + c p_{n,1} \mbox{,\ when\ } k =0.
\end{equation}

In section \ref{linear} we present a simple combinatorial
interpretation for this recurrence relation. Here are the
first five rows of $\Pi^{-1}$.
{\footnotesize\begin{equation*}
\begin{pmatrix}
1                & 0                 & 0        & 0 & 0 \\ 
c                & 1                 & 0        & 0 & 0 \\
c + c^2          & 1 + 2\,c          & 1        & 0 & 0 \\
c + 3\,c^2 + c^3 & 1 + 5\,c + 3\,c^2 & 2 + 3\,c & 1 & 0 \\
c + 6\,c^2 + 6\,c^3 + c^4 & 1 + 9\,c + 14\,c^2 + 4\,c^3 & 
3 + 11\,c + 6\,c^2 & 3 + 4\,c & 1 
\end{pmatrix}
\end{equation*}}

We can also obtain immediately the generating functions for
the sequences $\{p_{n,k}\}_n$ as in Haagerup-Thorbj{\o}rnsen
\cite{ht2}, \S 6. Let $P_k(z) = \sum_{n \geq 0} p_{n,k}
z^n$. Then $P_0(z) = \frac{1 - (c -1) z -\sqrt{(1 -az)(1 -
bz)}}{2 z}$ is the moment generating function of the
Marchenko-Pastur distribution. From equation (\ref{Prec1})
we get that for $k \geq 1$
\begin{equation}\label{Prec3}
c P_{k+1}(z) = (z^{-1} - (1 + c)) P_k(z) - P_{k-1}(z)
\end{equation}
and from equation (\ref{Prec2}) we get that
\begin{equation*}
\bigg(\frac{P_0(z) -1}{z}\bigg) = c P_0(z) + c P_1(z)
\end{equation*}
From this equation and the functional equation for $P_0$ we
get that \[ P_1(z) = \bigg( \frac{ P_0(z) - 1}{c}\bigg)
P_0(z) \] Then by induction we get, as noted in
\cite{ht1}, Lemma 6.3, from (\ref{Prec3}) that for all
$k$, \[ P_k(z) = \bigg( \frac{P_0(z) - 1}{c}\bigg)^k P_0(z) \]
In Remark \ref{lineargenremark} we give a simple and elegant
diagrammatic interpretation of this formula.

Finally let us note that the relation $2 T_n(x) = U_n(x) -
U_{n-2}(x)$ between the Chebyshev polynomials of the first
and second kind implies that
$$\label{gamma_pi}
\wt\Gamma_n(x) +
\wt\Gamma_{n-1}(x) = \Pi_n(x) - c \Pi_{n-2}(x)
$$
for $n \geq 2$ and thus (since $d_n + d_{n-1} = 0$ for $n \geq 3$)
\begin{equation}\label{firstsecond}
\Gamma_n(x) + \Gamma_{n-1}(x) = \Pi_n(x) - c
\Pi_{n-2}(x).
\end{equation}

\section{Wishart Matrices}\label{wishart}

Let $(\Omega, \Sigma, {\sf P})$ be a probability space and
$G : \Omega \rightarrow M_{M,N}(\C)$ be a random matrix with
entries $\{ g_{i,j} \}$.  Suppose that $\{ g_{i,j} \}$ are
independent complex Gaussian random variables with
$\E(g_{i,j}) = 0$ and $\E(|g_{i,j}|^2) = 1/N$ for all $1
\leq i \leq M$ and $1 \leq j \leq N$. Then $X = G^\ast G$ is
a particular case of a complex Wishart matrix. To simplify
the terminology we shall henceforth say that $X$ is a {\em
Wishart matrix} if $X = G^\ast G$ and $G$ is as above for
some $M$ and $N$.

We are interested in the behaviour of the eigenvalue
distribution of $X_{N}$ as $M$ and $N$ tend to
infinity. Thus we assume that we have sequences $\{ M_k
\}_k$ and $\{ N_k \}_k$ of positive integers such that $c :=
\lim_k M_k/N_k$ exists and $0 < c < \infty$. In order to get
a second order limiting distribution we shall further assume
that the limit $c' = \lim_{k \rightarrow \infty} M_k - c
N_k$ exists. Whenever asymptotics are discussed in this
paper we shall always assume that $M$ and $N$ are chosen
from sequences $\{M_k\}_k$ and $\{N_k\}_k$ satisfying the
two limiting behaviours above. When we take a limit as $k$
tends to infinity we shall denote this as $\lim_N$. For
further details and references see \cite{ht1}, \cite{ht2},
and \cite{mn}.

It was shown in Cabanal-Duvillard \cite{thierry} that the
family of random variables $\{ \Tr(\Gamma_n(X_{N})) \}_{n
> 0}$ is asymptotically Gaussian and independent as $N$
tends to infinity, where the $\Gamma_k$'s are the shifted Chebyshev
polynomials of the first kind constructed in section
\ref{firstkind}. We wish to extend this to a collection of
independent Wishart matrices. So suppose that for each $M$ and $N$
we have $G_1, G_2, \dots , G_p$ each as above but in addition such
that the entries of all the $G_i$'s are independent. Thus for each
$M$ and $N$ we have Wishart matrices $X_1^{(N)}, \dots , X_p^{(N)}$
where $X_i^{(N)} = G_i^\ast G_i$. We now wish to construct for each
$k$-tuple of positive integers $\vec m = (m_1, \dots , m_k)$ and
each sequence of indices $\vec i = (i_1, \dots , i_k)$ with $1 \leq
i_j \leq p$ and $i_j \not= i_{j+1}$ a random variable $S_{\vec m,
\vec i}$: $$ 
S_{\vec m, \vec i}^{(N)} = \Tr\big( \Pi_{m_1}(X_{i_1}^{(N)})
\Pi_{m_2}(X_{i_2}^{(N)}) \cdots \Pi_{m_k}(X_{i_k}^{(N)}) \big)
$$ 
where the $\Pi_k$'s are the shifted Chebyshev polynomials of the
second kind constructed in section \ref{pi-poly}. $S_{\vec m, \vec
i}^{(N)}$ only depends on the equivalence class of $(\vec m, \vec
\imath\,)$ under cyclic permutation, so in Theorem \ref{mt} below we
shall assume that we have chosen one representative from each
equivalence class. Let $|(\vec m, \vec \imath\,)|$ be the number of
cyclic equivalence classes, i.e. the number of $1
\leq l
\leq k$ such that for $1 \leq r \leq k$ we have $m_r = m_{k+l + r}$
and $i_r = i_{k+l + r}$, where the indices are taken modulo $k$.

\begin{theorem}\label{mt}
Let $\{\Gamma_n\}_n$ be the shifted Chebyshev polynomials of
the first kind as considered by Cabanal-Duvillard and let
$\{ \Pi_n \}_n$ be the orthogonal polynomials of the
Marchenko-Pastur distribution (which are shifted Chebyshev
polynomials of the second kind).  Let $X_1, \dots , X_p$ be
independent Wishart matrices and consider in addition to
$\Tr(\Gamma_n(X_i))$ also, for $k\geq 2$, the collection of
random variables $\Tr( \Pi_{m_1}( X_{i_1}) \cdots \ab
\Pi_{m_k}(X_{i_k}) )$, where the Wishart matrices which appear
must be cyclically alternating, i.e.,
$i_1\not=i_2\not=i_3\not=\dots\not=i_k \not=i_1$.  These
latter traces depend only on the equivalence class of $(i_1,
\dots , i_k)$ and $( m_1, \dots m_k)$ under cyclic
permutation. Assuming that we have chosen one representative
from each equivalence class, the random variables
$$\{ \Tr( \Gamma_n( X_i ) ) \} \cup \{ \Tr(
\Pi_{m_1}( X_{i_1}) \cdots \Pi_{m_k}(X_{i_k}) ) \}$$ are
asymptotically independent and Gaussian.

Moreover the limiting means of $\Tr (\Gamma_n(X_N) ) $ and
$S_{\vec m, \vec \imath}$ are given by
\[ \lim_N \E( \Tr (\Gamma_n(X_{N,i}) ) ) )=  (-1)^n c'
\mbox{ and } \]
\[
\lim_N \E( S_{\vec m, \vec \imath} ) = 0
\]

Finally the asymptotic variances of $S_{\vec m, \vec i}$ and 
$\Tr (\Gamma_n(X_N) ) $ are given by
\[
\lim_N  \kappa_2( \Tr (\Gamma_n(X_{N,i}), \Tr (\Gamma_n(X_{N,i}) )\,
) = n c^n \mbox{ and } \]
\[
\lim_N \kappa_2(Tr (S_{\vec m, \vec i}), Tr (S_{\vec m, \vec i})) =
\lim_N  \E( |\Tr (S_{\vec m, \vec i})|^2\, ) = |(\vec m, \vec
\imath\, )| c^{m_1 + \cdots + m_k}
\]
\end{theorem}

\begin{remark}
In proving Theorem \ref{mt} we shall show along the way that
\[
\lim_N \Tr (\Pi_n(X_{N,i}) ) ) = 
\begin{cases} 0 & n \mbox{ is even }\\
c' c^k & n = 2 k + 1 \\
\end{cases}
\]

\end{remark}

In Mingo and Nica \cite{mn}, \S 9 it was shown that the
traces of words in $\{X_i^{(N)}\}_i$ are asymptotically
Gaussian, in fact that all the cumulants of order higher than 2 are
asymptotically 0, thus the random variables
$\{\Tr(\Gamma_m(X_i^{(N)}))\}_{m,i}$ and $\{ S_{\vec m, \vec
i}^{(N)}\}_{\vec m, \vec i}$ are asymptotically Gaussian.

We shall prove Theorem \ref{mt} by showing that
asymptotically $ \{ S^{(N)}_{\vec m, \vec i} \}_{\vec m,
\vec i} \ab\cup\{ \Tr(\Gamma_m(X^{(N)}_j)) \}_{m,i} $ are
independent, i.e.  that they diagonalize the covariance
$\kappa_2$. In \cite{mn}, \S 7 it was shown that that
these covariances can be expressed in terms of planar diagrams
called non-crossing annular permutations. In the next two
sections we will show how to relate these diagrams to the
polynomials $\{ \Gamma_n \}_n$ and $\{ \Pi_n \}_n$.

\section{Non-Crossing Circular Half-Permutations}\label{circular}

We introduce the notion of a non-crossing circular
half-permutation.  A non-crossing circular half-permutation
on $[m]$ is a non-crossing permutation $\pi$ of $[m]$ together
with a selection of blocks of $\pi$ that satisfy a
condition described below.

We shall begin by recalling the definition of non-crossing
permutations. Let $[m] = \{1, 2, 3, \dots , m\}$. Let $\pi$
be a partition of $[m]$. We say $\pi$ has a crossing if
there are $r < s < t < u$ with $r$ and $t$ belonging to one
block of $\pi$ and $s$ and $u$ belonging to another. We say
$\pi$ is {\it non-crossing} if it has no crossings.

Another useful picture is to consider $\pi$ as a permutation
of $[m]$. Each block of $\pi$ is arranged into increasing
order and these form the cycles of a permutation of
$[m]$. The permutations so obtained are characterized by the
equality $\#(\pi) + \#(\gamma_m\, \pi^{-1}) = m+1$, where
$\gamma_m$ is the permutation with one cycle $(1, 2, 3,
\dots , m)$ and $\#(\sigma)$ is the number of cycles in the
permutation $\sigma$; see Biane \cite{biane}. As there is a
bijection between permutations of $[m]$ satisfying this
geodesic condition and partitions satisfying the
non-crossing condition we shall denote them both by
$\pi$. When it is necessary to emphasize the distinction we
shall denote the non-crossing partitions by $NC(m)$ and the
non-crossing permutations by $S_{NC}(m)$.

Perhaps the simplest description however is in terms of
planar diagrams. Given a partition $\pi$ we place the
numbers 1, 2, 3, \dots , $m$ around a circle in clockwise
order, and in the interior of the circle connect the points in the
same block. If this can be drawn so that the blocks do not cross
then the partition is non-crossing.

Given a non-crossing partition $\pi$ there is another
partition called the Kreweras complement which we shall
denote by $\pi^c$.  The complement can be described several
ways. First, let us regard $\pi$ as a partition of
$[m]$. Let us consider another set $[\ol m] = \{\ol 1, \ol
2, \dots , \ol m\}$ and arrange the union $[\ol m] \cup [m]$
in the order $\ol 1 < 1 < \ol 2 < 2 < \cdots < \ol m <
m$. Then $\pi^c$ is the largest partition of $[\ol m]$ such
that $\pi^c \cup \pi$ is a non-crossing partition of $[\ol
m] \cup [m]$.

Alternatively we can regard $\pi$ as a permutation of $[m]$
and then $\pi^c = \gamma_m \pi^{-1}$, see Biane
\cite{biane}.

\begin{definition}\label{circulardef}
A {\it non-crossing circular half-permutation} of $[m]$ is a
tuple $(\pi, \ol B, B_1, B_2, \dots , B_k)$ where $\pi$ is a
non-crossing permutation of $[m]$, $ k \geq 1$, $B_1$, \dots $B_k$
are blocks of $\pi$, and $\ol B$ is a block of $\pi^c$ such that
$B_i \cap \bar B \not= \emptyset$ for $i = 1, \dots ,
k$. The blocks $B_1$ , \dots , $B_k$ are called the {\it
open blocks} of $\pi$ and the remaining blocks of $\pi$ are
called the {\it closed blocks}. Frequently and when no confusion
can occur we shall simply denote $(\pi, \ol B, B_1, \dots ,
B_k)$ by $\pi$. The set of non-crossing circular
half-permutations with $k$ open blocks will be denoted 
$\ncc(m)_k$.\end{definition}

\begin{lemma}\label{blocksize}
Suppose that $\pi \in NC(m)$, $B$ is a block of $\pi$, $\ol
B$ is a block of $\pi^c$ such that $B \cap \ol B \not=
\emptyset$.  Then $|B \cap \ol B| = 1$.
\end{lemma}

\begin{proof}
Suppose $i, j \in B \cap \ol B$ with $i < j$. For this proof
we shall use the picture of $\pi^c$ as a partition on $[\ol
m]$ such that $\pi^c \cup \pi$ is a non-crossing partition
of $[\ol m] \cup [m]$. Then $\ol i < i < \ol j < j$ and $i,j
\in B$ and $\ol i, \ol j \in \ol B$ which gives a crossing
of $\pi^c \cup \pi$. Hence $|B \cap \ol B| = 1$.
\end{proof}

\begin{definition}\label{initial}
Let $\pi \in NC(m)$, $B \in \pi$, and $\ol B \in \pi^c$ with
$B \cap \ol B = \{ i \}$. Then $i$ is called the {\it
initial point} of $B$ {\it relative} to $\ol B$. We call
$\pi^{-1}(i)$ the {\it final point} of $B$ {\it relative} to
$\ol B$. If $B$ is a singleton then its initial and final
points will coincide.
\end{definition}

We can now present another picture of non-crossing circular
half-permutations that explains the terminology we have
introduced above.

In \cite{mn} it was shown that for a sequence of Wishart
matrices $\{ X_{N} \}$ the correlation of the moments
$$
\E(\Tr(X^m_N)\, \Tr(X^n_N)) - \E(\Tr(X^m_N))\, \E(\Tr(X^n_N))
$$ 
converged as $N$ tends to infinity, to $\ds \sum_{\pi \in
S_{NC}(m,n)} c^{\#(\pi)}$. $S_{NC}(m,n)$ is the collection
of non-crossing annular permutations on the $(m,n)$-annulus
and $\#(\pi)$ is the number of cycles of $\pi$.  These were
the subject of \cite{mn}; but we shall recall the pertinent
facts here.

The notion of a non-crossing annular partition extends to
the annulus the idea of a non-crossing partition on a
disc. Given two integers $m$ and $n$ and two concentric
circles we have on the outer one the points 1, 2, \dots ,
$m$ in clockwise order and on the inner circle we have the
points $m+1$, \dots , $m + n$ in counter-clockwise order. We
call this the $(m,n)$-{\it annulus}. We want to study
partitions of $[m + n]$ such that when drawn on the
$(m,n)$-annulus there is at least one block connecting the
two circles and the blocks do not cross.  Being non-crossing
on the $(m,n)$-annulus is weaker than being non-crossing on the
$m+n$-disc.

Given a non-crossing partition of the disc we can always put
the points of each block in standard order and obtain a
permutation satisfying Biane's condition. In the case of the
annulus we have to distinguish between permutations whose
orbits as a set are the same but which visit the points in different
orders.

Informally a permutation $\pi$ in $S_{m+n}$ is non-crossing
on the $(m,n)$-annulus if we can connect the points in
cycles in the order visited by $\pi$ in such a way that the
blocks do not cross or self-intersect and enclose their
interior in the clockwise orientation..

\noindent
$\vcenter{\hsize151pt\includegraphics{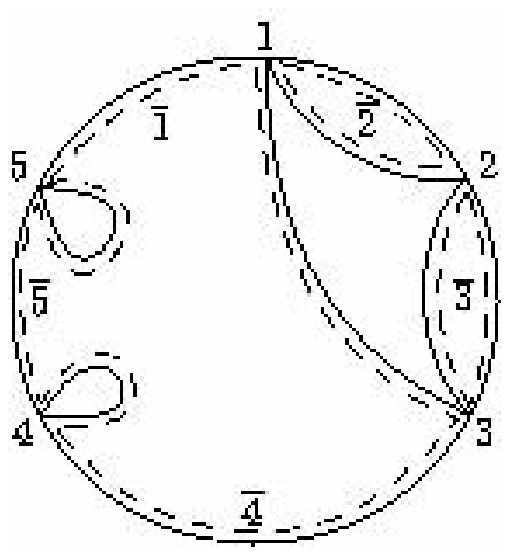}}$\hfill
$\vcenter{\hsize 180pt\raggedright\parindent0pt\small\vfill
{\bf Figure \figno.} The non-crossing partition $\pi = (1, 2,
3)(4)(5)$ and its complement $\pi^c$ = (1, 4, 5) (2) (3).
\vfill}$
%\figbox{{\small $\pi = (1,2,3)(4,6),\ab (5),\ab (7,8)$}}

Elements of $S_{NC}(m,n)$ are permutations of $[m+n]$ which
satisfy an annular geodesic condition similar to the one
found by Biane. Let $\gamma_{m,n}$ be the permutation of
$[m+n]$ with two cycles 
$$
\gamma_{m,n} = (1,2,3, \dots , m)(m+1, \dots , m+n)
$$ 
A permutation $\pi$ is {\it non-crossing annular} if it connects the
two circles and 
\begin{equation}\label{geodesic}
\#(\pi) + \#( \gamma_{m,n} \pi^{-1}) = m + n
\end{equation}
There are a variety of equivalent ways of expressing
this condition; see \cite{mn}, \S 3, 4, 5, 6 and
\cite{ms}, \S 2.

\noindent
$\vcenter{\hsize175pt\includegraphics{fig2.eps}}$
\hfill
$\vcenter{\hsize 150pt\raggedright\parindent0pt\small\vfill
{\bf Figure \figno.} $\pi$ = (1,2,3,12) (4,9) (5,6,7) (8) (10,
11) is a non-crossing annular (8,4)-permutation and its
complements $\pi^c = (1, 9, 5, 8)\ab (2) (3) (4,10,12) (6)
(7)$. \vfill}$

We now wish to describe what happens if we cut one of these
annular diagrams into two pieces which we shall call {\it
half-permutations}. Informally we mean we separate the two
circles; blocks that lie completely on one circle are called
{\it closed blocks} and blocks that connect the two circles
are cut in half and are called {\it open blocks}. Let us
consider the example $\pi = (1,2,3,12) \, (4,9)\,
(5,6,7)\,\ab (8) (10,11)$ which is a non-crossing permutation of
the (8,4)-annulus. There are two {\it through-blocks}
i.e. blocks that connect the two circles: (1,2,3,12) and
(4,9). So if we consider the half on the outside circle we
have the half-permutation
$\{(1,2,3)\,(4)\,(5,6,7)\,(8)\}$. The open blocks are
(1,2,3) and (4) and the closed blocks are (5,6,7) and (8).

\setbox1=\vbox{\hsize300pt\parindent0pt\small
{\bf Figure \figno.} The two non-crossing circular
half-permutations arising from $\pi = (1,2,3,12) (4,9) (5,6,7)
(8) (10, 11)$.}
\begin{center}
\includegraphics{fig4.eps}\qquad\qquad
\includegraphics{fig5.eps}

\

\leavevmode\box1
\end{center}

We now give a formal definition.

\begin{notation}
Let $S_{NC}(m,n)$ be the set of non-crossing annular 
permutations. For $\pi \in S_{NC}(m,n)$ we call
$\gamma_{m,n} \pi^{-1}$ the {\it annular complement} of
$\pi$ and denote it also by $\pi^c$. $\pi^c$ is also non-crossing
by  equation (\ref{geodesic}). For $\pi \in
S_{NC}(m,n)$ and $B$ a cycle of $\pi$ we say that $B$ is
{\it through-block} if $B \cap [m] \not= \emptyset$ and $B
\cap [m+1, m+n] \not = \emptyset$
\end{notation}

\begin{notation}
Let $\pi$ be a permutation of $[n]$ and let $A \subset
[n]$. Let $\pi_{|_A}$ denote the permutation induced by
$\pi$ on $A$, i.e. for $i \in A$, $\pi_{|_A}(i) = j$ if
$\pi^k(i) = j \in A$ and $\pi^l(i) \not\in A$ for $1 \leq l
< k$.
\end{notation}

\begin{lemma}\label{lemma5}
Let $\pi \in S_{NC}(m, n)$ and $\pi_1 = \pi_{|_{[m]}}$. Let
$\ol B_1, \ol B_2, \dots , \ol B_k$ be the through-blocks of
$\pi^c$. Let $\ol B = (\ol B_1 \cup \cdots \cup \ol B_k ) \cap [m]$,
then
$\ol B$ is a cycle of $\pi_1^c$.
\end{lemma}

\begin{proof}
We must show that $\ol B$ is invariant under $\pi_1^c$ and that
$\pi_1^c$ acts transitively on $\ol B$. Let $i \in \ol B$. First
suppose $\pi^c(i) \in [m]$, then $\pi^{-1}(i) \in [m]$ and
so $\pi_1^{-1}(i) = \pi^{-1}(i)$; thus $\pi_1^c(i) =
\gamma_m \pi_1^{-1}(i) = \pi^c(i) \in \ol B$.

Suppose $\pi^c(i) \in [m+1, m+n]$, then $\pi^{-1}(i) \in [m
+ 1, m + n]$.  Let $k$ be such that $\pi^{-k}(i) \in [m]$
but $\pi^{-l}(i) \in [m + 1, m + n]$ for $1 \leq l < k$. Let
$j = \pi^{-k+1}(i)$.Then $\gamma_{m,n} \pi^{-1} (j) =
\gamma_m \pi_1^{-1} (i)$. So $\pi_1^c(i)$ is in a
through-block of $\pi^c$ and thus $\pi_1^c(i) \in \ol B$. Hence
$\ol B$ is invariant under $\pi_1^c$.

If $\ol B_i$ is a through-block of $\pi^c$ and $i \in \ol B_i$ is
such that $(\pi^c)^{-1}(i) \in [m + 1, m + n]$ then we shall
call $i$ the {\it initial point} of $\ol B_i$; conversely if
$\pi^c(i) \in [m + 1, m + n]$ then we shall call $i$ the
{\it final point} of $\ol B_i$.  We have just seen in the first
paragraph of the proof that if $i \in \ol B_i$ is not the final
point then $\pi_1^c$ moves through the points of $\ol B_i \cap
[m]$ until it comes to the final point of $\ol B_i$. In the
second paragraph we saw that $\pi_1^c$ takes the final point
of one block of $\pi^c$ to the initial point of another
block.

Each of the blocks $\ol B_i$, being a through-block, has a
non-empty intersection with $[m + n]$. By property (AC-3) of
\cite{mn}, Def. 3.5, applied to $\pi^c$ we see that
once $\gamma_m$ visits a block $\ol B_i$ it doesn't visit it
again until all the other blocks $\ol B_j$ ($j \not= i$) have
been visited. So let us order the blocks $\ol B_i$ so that
$\gamma_m$ visits the blocks in the order $\ol B_1$, $\ol
B_2$, \dots , $\ol B_k$. If $i \in \ol B_r \cap [m]$ and $i$
is a final point of $\ol B_i$ then $i$ is an initial point of
a block $C$ of $\pi$ and so $\pi_1^{-1}(i)$ is the final point
of $C$; thus $\gamma_m \pi_1^{-1}(i)$ is the initial point of
$\ol B_{r+1}$. \end{proof}

\begin{definition}
Let $\pi$ be a non-crossing annular $(m,n)$-permutation. We shall
construct a pair $\pi_1$ and $\pi_2$ of non-crossing circular
half-permutations from $\pi$.

Let $\pi_1 = \pi_{|_{[m]}}$. By
\cite{mn}, Def.  3.5, (AC-1) $\pi_1$ is non-crossing.

Suppose the through-blocks of $\pi^c$ are $\ol B_1$, $\ol
B_2$, \dots , $\ol B_r$. Let $\ol B = (\ol B_1 \cup \cdots
\cup \ol B_r) \cap [m]$. Let $B_1, B_2, \dots , B_s$ be the through-blocks of
$\pi$. By Lemma \ref{lemma5}, $\ol B$ is a
cycle of $\pi_1^c$ and the blocks $B_i \cap [m]$ all meet $\ol
B$. Thus $(\pi_1, B_1 \cap [m], B_2 \cap [m], \dots , B_s
\cap [m], \ol B)$ is a non-crossing circular
half-permutation.

Similarly let $\pi_2 = \pi_{|_{[m + 1, m + n]}}$. Then $\pi_2$ is
non-crossing on $[m + 1, m + n]$. Let $\ol B' = (\ol B_1
\cup \cdots \cup \ol B_r)\cap [m + 1, m + n]$. Then $(\pi_2,
B_1 \cap [m + 1, m + n], B_2 \cap [m + 1, m + n], \dots ,
B_s \cap [m + 1, m + n], \ol B')$ is a non-crossing circular
half-permutation on $[m + 1, m + n]$. $\pi_1$ and $\pi_2$
are the two half-permutations obtained from $\pi$.
\end{definition}

\begin{definition}
Let $(\sigma, B_1, B_2, \dots , B_k, \ol B)$ be a non-crossing
circular half-permutation of $[m]$. We shall say that the cycles
$B_1, B_2, \dots , B_k$ are in {\it cyclic order} if for each $i$
and each $x \in B_i$ and $r$ such that $\gamma_m^r(x) \in \cup_{l
\not = i} B_l$ and $\gamma_m^s(x) \not \in \cup_{l \not = i} B_l$ for
$1 \leq s < r$, $\gamma_m^r(x) \in B_{i+1}$ where we interpret
$B_{k+1}$ as $B_1$. 

\end{definition}

\begin{lemma}
Let $(\sigma, B_1, B_2, \dots , B_k, \ol B)$ be a non-crossing
circular half-permutation of $[m]$. There is an ordering of the
cycles $B_1, B_2, \dots, B_k$ such that they are in cyclic order.

\end{lemma}

\begin{proof}
We showed in Lemma \ref{blocksize} that for each $i$, $|B_i \cap \ol
B| = 1$. Let us label these points $x_1, x_2, \dots , x_k$, i.e.
$\{x_i \} = |B_i \cap \ol B|$. Moreover suppose that the
cycles $B_1, B_2, \dots , B_k$ are ordered so that $\gamma_m$ visits
the $x_i$'s in the order $x_1, x_2, \dots , x_k$. We shall show
that this puts the cycles into cyclic order.

Choose $i$ and $x \in
B_i$ and let $r$ be such that $\gamma_m^r(x) \in \cup_{l \not =
i}B_l$ and $\gamma_m^s(x) \not \in \cup_{l \not = i} B_l$ for $1
\leq s < r$. Since the points  $\ol x_1, x_1, \ol x_2, x_2, 
\dots , \ol x_k, x_k$, of $[\ol m] \cup [ m ]$, are in cyclic order
we can only have $x_i \leq x < x_{i + 1}$, otherwise there would be
a crossing between $B_i$ and $\ol B$. Hence $\gamma_m^r(x) =
x_{i+1} \in B_{i+1}$. 
\end{proof}

\begin{notation}
Let $(\sigma, B_1, B_2, \dots , B_k, \ol B)$ be a non-crossing
circular half-permutation of $[m]$. We denote $x_i$ the initial point
of $B_i$ relative to $\ol B$ and by $y_i$ the final point of $B_i$. 
\end{notation}

\begin{remark}
In the following theorem we have two non-crossing circular
half-permutations each with $k$ open blocks. We show that one can
construct $k$ non-crossing annular permutations from the pair of
half-permutations. The diagrammatic interpretation is very simple.
We fix an open block on the first half-permutation and then there
are $k$ open blocks on the second half-permutation to connect it to.
Once this choice has been made there are no further choices for
pairing up the remaining $k-1$ open blocks on each half-permutation.
\end{remark}

$$
\vcenter{\hsize75pt\includegraphics{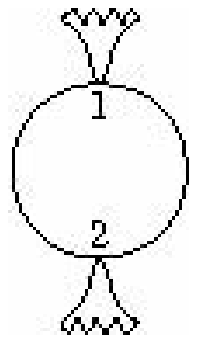}} \mbox{ and }
\vcenter{\hsize75pt\includegraphics{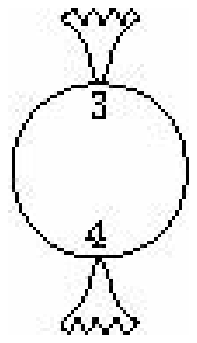}} \mbox{ produce }
$$

\medskip{\parindent0pt$%
\vcenter{\hsize155pt\includegraphics{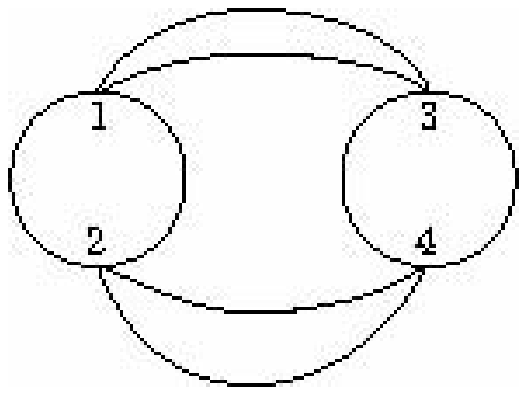}} \mbox{ and
}
\vcenter{\hsize150pt\includegraphics{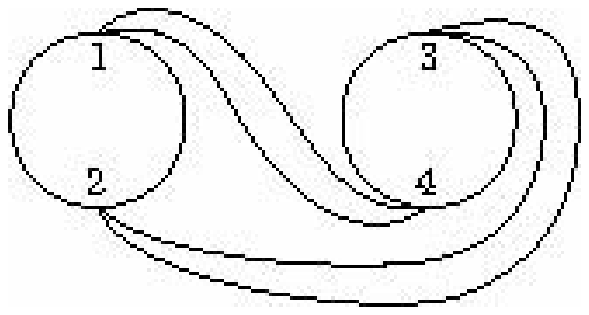}}
$}

\begin{theorem}
Let $(\sigma, B_1, \dots , B_k,\ol B)$ and $(\tau, C_1, \dots , C_k
, \ol C)$ be non\-cross\-ing circular half-permutations of $[m]$ and
$[n]$ respectively, each with $k$ open blocks.

Then there are
exactly $k$ non-crossing $(m,n)$-annular permutations $\pi^{(1)}$,
$\pi^{(2)}$, \dots , $\pi^{(k)}$ such that for each $i$,
$\pi^{(i)}_1 = \sigma$ and $\pi^{(i)}_2 = \tau$.

\end{theorem}

\begin{proof}

Let $1 \leq s \leq k$; we shall construct $\pi^{(s)}$. 
We shall regard $\tau$ as a non-crossing permutation of $[m + 1,
m + n]$ by identifying $[n]$ with its translate in $[m +1, m + n]$.
Moreover shall denote by $\tilde \gamma_n$ the permutation of $[
m+1, m+n]$ with the one cycle $(m+1, m+2, \dots , m+n)$. Let the open
blocks $B_1$, $B_2$. \dots , $B_k$ of $\sigma$ and the open blocks
$C_1$, $C_2$, \dots , $C_k$ of be arranged in cyclic order. Let
$x_1, x_2, \dots , x_k$ and $y_1, y_2, \dots , y_k$ be such that
$B_i \cap \ol B = \{x_i \}$ and $C_i \cap \ol C = \{ y_i \}$. Let  $$
\pi^{(s)} = (x_1, y_{k-1+s}) (x_2, y_{k-2+s}) \cdots
            (x_k , y_s) \sigma \tau
$$
where the index of $y_{k-i+s}$ is interpreted modulo $k$. Let us
show that $\#(\pi^{(s)}) + \#(\gamma_{m,n} (\pi^{(s)})^{-1}) = m +
n$. This will show that $\pi^{(s)}$ is a non-crossing annular
permutation by the geodesic condition, \cite{mn}, Theorem
6.1.

Let $\pi_1 = (x_1, y_{k-1+s}) \sigma \tau$. Then $\#(\pi_1) =
\#(\sigma) + \#(\tau) -1$ as $B_1$ and $C_{k-1+s}$ are disjoint
cycles of $\sigma \tau$. Also 
$$
\#(\gamma_{m,n} \pi_1^{-1}) = \#(\gamma_{m,n} \sigma^{-1} \tau^{-1}
(x_1 y_{k-1+s})) = \#(\gamma_m \sigma^{-1}) + \#(\tilde \gamma_n
\tau^{-1}) - 1
$$ 
as $\ol B$ and $\ol C$ are disjoint cycles of $\gamma_{m,n}
\sigma^{-1} \tau^{-1}$. Thus by the geodesic condition
$\pi_1$ is non-crossing annular.

Now let $\pi_r = (x_1, y_{k-1+s}) \cdots (x_r, y_{k-r+s})
\sigma \tau$. We shall show by induction that $\pi_r$ is
non-crossing annular. Since $B_i$ and $C_{k-i+s}$ are
disjoint cycles of $\pi_{i-1}$ we have that $\#(\pi_r) =
\#(\pi_{r-1}) -1$. We shall complete the proof by showing
that $\#(\gamma_{m,n} \pi_r^{-1}) = \#(\gamma_{m,n}
\pi_{r-1}^{-1}) + 1$, as this will inductively show that
$\pi_r$ satisfies the geodesic condition.

Write $\sigma = \sigma_0 \ol B$ and $\tau = \tau_0 \ol C$, then
$$
\gamma_{m,n} \pi_r^{-1} = \gamma_{m,n}\sigma_0^{-1} \tau_0^{-1}
(\ol{B})^{-1} (\ol C)^{-1} (x_1, y_{k-1+s}) \cdots (x_r, y_{k-r+s})
$$
Thus it remains to shown that $x_r$ and $y_{k-r+s}$ are in the same
cycle of 
$$
B^{-1} C^{-1} (x_1, y_{k-1+s}) \cdots (x_{r-1}, y_{k-r+1+s})
$$ 
or equivalently that $x_r$ and $x_{k-r+s}$ are in the same cycle of 
$$
(x_1, y_{k-1+s}) \cdots (x_{r-1}, y_{k-r+1+s})B C
$$ 
We establish this by recursively applying the following sublemma.

\medskip\noindent{\it
Sublemma}\/ Suppose $b_1 < b_2 < \cdots < b_k < c_1 < c_2 < \cdots <
c_l$; $1 \leq i_1 < i_2 < \cdots < i_t \leq k$; $l \geq j_1 > j_2
> \cdots > j_t \geq 1$; and $D = (b_1, b_2, \dots , b_k, c_1, c_2,
\dots , c_l)$ be a cycle. For $r \not = s$, $b_{i_s}$ and $c_{j_s}$
are in the same cycle of $(b_{i_r}, c_{j_r}) D$. 

\medskip\noindent{\it
Proof of sublemma}\/ By direct calculation, 
$$
(b_{i_r}, c_{j_r}) D =
(b_1, \dots , b_{i_{r-1}}, c_{j_r}, \dots , c_l) 
(b_{i_r}, \dots , b_k, c_1, \dots , c_{j_{r-1}})
$$
If $s < r$ then $b_{i_s}$ and $c_{j_s}$ are in the first cycle; if
$s > r$ they are in the second. %$\vartriangleright$

To complete the proof we must show that every non-crossing
annular permutation arises from a pair of non-crossing
circular half-permutations. So let $\pi$ be a non-crossing
annular permutation. We must show that when we cut $\pi$
into $\pi_1$ and $\pi_2$, a pair of non-crossing circular
half-permutations, and then reassemble them we recover
$\pi$. The only point that needs to be checked concerns the
through blocks of $\pi$. So suppose that $D_1$, $D_2$, \dots
, $D_k$ are the through blocks of $\pi$ and we let $B_1 =
D_1 \cap [m]$, $B_2 = D_2 \cap [m]$, \dots , $B_k = D_k \cap
[m]$ be the open blocks of $\pi_1$ and $C_1 = D_k \cap [m+1,
m+m]$, $C_2 = D_{k-1} \cap [m=1, m+n]$, \dots , $C_k = D_1
\cap [m+1, m+n]$. Let us further suppose that the through
blocks $D_1$, \dots , $D_k$ have been ordered so that $B_1$,
\dots , $B_k$ are in cyclic order. By \cite{mn}, Def. 3.5
(AC-3) (see the middle diagram of \cite{mn}, Figure
3.7) the blocks
$C_k$, \dots , $C_1$ are in cyclic order. Let $x_1$, \dots ,
$x_k$ be the initial points of $B_1$, \dots , $B_k$; and
similarly $y_1$, \dots , $y_k$ the initial points of $C_1$,
\dots $C_k$. Then by \cite{mn}, Def. 3.5 (ANS-2), $D_i =
(x_i, y_{k-i+1}) B_i C_{k-i+1}$ for $1 \leq i \leq k$. Thus
$\pi = (x_1, y_k) \cdots (x_k, y_1) \pi_1 \pi_2$.  \end{proof}

\setbox1=\vbox{\hsize300pt\parindent0pt\small
{\bf Figure \figno.} The two half-permutations of figure 9 
reassembled to give the original permutation.}

\bigskip\begin{center}
\includegraphics{fig3.eps}

\

\leavevmode\box1
\end{center}

\section{Non-crossing Linear Half permutations}\label{linear}

In the previous section we saw that a non-crossing annular
permutation could be decomposed into two non-crossing
circular half-permutations.  We shall need to decompose
these circular half-permu\-ta\-tions yet further into what we
shall call non-crossing linear half-permutations.

Suppose that $m = m_1 + m_2 + \dots m_k$, with each $m_i >
0$; and that we have intervals $I_1, \dots , I_k \subset
[m]$ of cyclically consecutive points. Moreover suppose that
the interval $I_j$ has $m_j$ points each of colour $i_j$ and
that cyclically adjacent intervals have different colours.

Let $(\pi, B_1, \dots , B_k, \ol B)$ be a non-crossing
circular half-permutation of $[m]$ such that all points of a
cycle have the same colour and that each interval $I_r$
meets at least one open block of $\pi$. Then each cycle of
$\pi$ can meet only one interval and so for each $r$,
$\pi|_{I_r}$ is a non-crossing permutation of
$I_r$. Moreover $\ol B$ meets the initial point of $I_r$. We
shall formalize this notion in the definition below.

\begin{definition}
A {\it non-crossing linear half-permutation} of $[m]$ is a
non-crossing circular half-permutation $(\pi , B_1,
\dots , B_k, \ol B)$ in which $1 \in \ol B$.  When $k = 0$ we
understand this to mean a non-crossing partition. We will denote by
$\ncl(m)_k$ the set of non-crossing linear half-permutations on
$[m]$ with $k$ open blocks.

If $I$ is a finite interval, a non-crossing linear
half-permutation of $I$ is a non-crossing circular
half-permutation $(\pi, B_1, \dots , B_k, \ol B)$ in which the
initial point of $I$ is in $\ol B$.

\end{definition}

We summarize our discussion  in the following theorem.

\begin{theorem}\label{decomposition}
Suppose $m  = m_1 + \cdots + m_k$ and
we have intervals $I_1, \dots , I_k$ of cyclically adjacent
points with the interval $I_r$ of length $m_r$ and coloured $i_r$
and the colours $i_1,  \dots , i_k$ cyclically alternating. Let
$(\pi, B_1, \dots B_k, \ol B)$ be a non-crossing circular
half-permutation of $[m]$ where each interval $I_r$ meets one
of the open blocks. Then each of $\pi|_{I_r}$ is a non-crossing
linear half-permutation whose open blocks are those of $\pi$
which meet $I_r$.
\end{theorem}

Given $0 \leq k \leq n$, let 
$$\ol p_{n,k} = \sum_{\pi \in \ncl(n)_k}
c^{\#(\pi)_{\rm cl}}$$ 
where $\#(\pi)_{\rm cl}$ denotes the number of closed blocks of
$\pi$. In this section we shall show that 
$$
\ol p_{n+1,k} = \ol p_{n, k-1} + (1 + c)\, \ol p_{n,k} +
c \ol p_{n, k+1}
$$
This is the same recursion as for the $p_{n,k}$'s (see
equation (\ref{Prec1})). 

\begin{theorem}\label{linearrecursion}
For $0 < k \leq n$
$$
\ol p_{n+1,k} = \ol p_{n, k-1} + (1 + c)\, \ol p_{n,k} +
c \ol p_{n, k+1}
$$
and for $k = 0$
$$
\ol p_{n+1,0} = c \ol p_{n,0} + c \ol p_{n,1}
$$
Moreover $\ol p_{n,k} = p_{n,k}$ for $0 \leq k \leq n$. 
\end{theorem}

\begin{proof}
We must show
that we can write $\ncl(n+1)_k$ as a disjoint union of the four
subsets  $$\ncl(n+1)_{k,1},\ \ncl(n+1)_{k,2},\ \ncl(n+1)_{k,3},\
\ncl(n+1)_{k,4}$$ and exhibit bijections \begin{eqnarray*}
&&\phi_1: \ncl(n+1)_{k,1} \rightarrow \ncl(n)_{k-1}\\
&&\phi_2: \ncl(n+1)_{k,2} \rightarrow \ncl(n)_{k}\\
&&\phi_3: \ncl(n+1)_{k,3} \rightarrow \ncl(n)_{k}\\
&&\phi_4: \ncl(n+1)_{k,4} \rightarrow \ncl(n)_{k+1}\\
\end{eqnarray*}
such that
\begin{eqnarray*}
\#( \phi_i(\pi))_{\rm cl} &=& \#(\pi)_{\rm cl} \mbox{\ for\ } i=
1,2 \mbox{ and } \\
1 + \#( \phi_i(\pi))_{\rm cl} &=& \#(\pi)_{\rm cl} \mbox{\ for\ }
i= 3,4
\end{eqnarray*}

Given $\pi$ in $\ncl(n+1)_k$ we look at the block containing
$n + 1$. There are four possibilities
\begin{itemize}
\item[-\ ] $n + 1$ is in an open block which is a singleton;

\item[-\ ] $n + 1$ is in an open block which is not a singleton;

\item[-\ ] $n+1$ is in a closed block which is a singleton; and

\item[-\ ] $n + 1$ is in a closed block which is not a singleton.
\end{itemize}

These four subsets clearly partition $\ncl(n+1)_k$. Next
we describe the maps $\{\phi_1, \phi_2, \phi_3, \phi_4\}$. In
cases (1) and (3), $n + 1$ is a singleton and $\phi_1$ and
$\phi_3$ remove this singleton leaving the other blocks alone.
In case (1) the number of open blocks decreases by one and
the number of closed blocks is constant; so
$\#(\phi_1(\pi))_{\rm cl} = \#(\pi)_{\rm cl}$. In case (3) the
number of open blocks is unchanged but the number of closed
blocks decreases by one so $1 +  \#(\phi_3(\pi))_{\rm cl} =
\#(\pi)_{\rm cl}$. In case (2) $n + 1$ is part of an open block,
we remove $n+1$ leaving the block open; thus both the number
open and the number of closed blocks is unchanged. 

For the case of (4) $n + 1$ is part of a closed block. We
remove $n+1$ and make the block open. Since $n+1$ is at an
endpoint it cannot be covered by another block and thus we
will not create a crossing by opening this block. However
the number of open blocks increases by one and the number of
closed blocks decreases by one and so $1 +
\#(\phi_4(\pi))_{\rm cl} = \#(\pi)_{\rm cl}$.

It is easy to see that the maps in (1), (2), and (3) are
bijections.  In the case of (4) we just have to show that
the construction can be reversed. Given a non-crossing
linear half-permutation on $[n]$ with $k+1$ open blocks we
add a point at $n+1$, make it part of the rightmost open
block, and then make this block closed.

In the case of $k = 0$ we just have two cases (3) and (4)
and we use the same maps as above.

We have now established that $\{ \ol p_{n,k} \}_{n,k}$
satisfy the same recursion as $\{ \ol p_{n,k} \}_{n,k}$, to
prove that they are equal we have to show that they are
equal for $n = 1$.

By direct calculation we have $p_{1,0} = c$ and $p_{1,1} =
1$. On the other hand when $n = 1$ a partition can have only
one block so $\ol p_{1,0} = c$ and $\ol p_{1,1} = 1$.
\end{proof}

We shall use the following theorem to calculate the limiting
mean of $\E(\Tr( \Gamma_n(X_N)))$ and $\E(\Tr(\Pi_n(X_N)))$
as $N \rightarrow \infty$.

\begin{theorem}\label{lineardecomp}
$$
\sum_{k=0}^{ [\frac{n-1}{2}] } c^k p_{n, 2k + 1} =
\sum_{\pi \in NC(n) } \#(\pi)\, c^{ \#(\pi) - 1}
$$
\end{theorem}

\begin{proof}
Let $\wt {\it NC}(n)$ be the set of pairs $( \pi , B)$ where
$\pi \in {\it NC}(n)$ and $B$ is a block of $\pi$. Let
$\psi(\pi, B) = c^{\#(\pi) - 1}$. Then the right hand side
can be written $\sum_{(\pi, B) \in \wt {\it NC}(n) }
\psi(\pi, B)$.

Let us define a map
\[
{\ds \bigcup_{k=0}^{ [\frac{n-1}{2} ] }} \ncl(n)_{2k+1}
\rightarrow \wt {\it NC}(n)
\]
which we shall denote $(\pi, B_1, \dots , B_{2k+1}) \mapsto
(\ol\pi, B)$ as follows. Arrange the open blocks $B_1,
\dots, B_{2k+1}$ from left to right as follows. $\min(B_1) <
\min(B_i)$ for $i > 1$ and $\max(B_i) < \max(B_{2k+1})$ for
$i < 2k+1$. Then $\min(B_2) < \min(B_i)$ for $i > 2$ and
$\max(B_i) < \max(B_{2k})$ for $i < 2k$. Then we join $B_1$
with $B_{2k+1}$, $B_2$ with $B_{2k}$, \dots, and $B_{k-1}$
with $B_{k+1}$. Call this new partition $\ol \pi$ and let $B
= B_k$. As we joined these blocks from the outside in we
have that $\ol \pi$ is non-crossing. Thus $(\ol \pi, B) \in
\wt{\it NC}(n)$.

Conversely starting with $(\ol \pi, B) \in \wt {\it NC}(n)$,
let $k$ be the number of blocks covering $B$. We split each
of these $k$ into two and declare them open. The open blocks
of $\pi$ will be $B$ and these $2k$ blocks.

If $(\pi, B_1, \dots , B_{2k+1}) \in \ncl(n)_{2k+1}$ with
$j$ closed blocks then $k + \#(\pi)_{\cl} = k +j$ and
$\#(\ol \pi) -1 = \#(\pi) -k - 1 = k + j$, so $c^k\,
c^{\#(\pi)_{\cl} } = \psi(\ol \pi, B)$. Summing over
$\cup_{k=0}^{ [\frac{n-1}{2} ] } \ncl(n)_{2k+1}$ now gives
us the result.
\end{proof}

\begin{remark}\label{lineargenremark}
The equation 
$$
P_k(z) = \Bigg(\frac{P_0(z) - 1}{c}\Bigg)^k P_0(z)
$$ 
has the following combinatorial interpretation. Every
non-crossing linear half-permutation with $k$ open blocks
can be written as the concatenation of $k$ non-crossing
linear half-permutations with one open block and one
non-crossing linear half-permutation with zero open blocks.
Indeed, given $\pi$ a non-crossing linear half-permutation
with $k$ open blocks, let $\pi_1$ be all the blocks of $\pi$
starting from the left up to and including the first open
block. Then $\pi_2$ starts after this open block and
continues up to the second block, and so on until all the
open blocks are exhausted; what remains is a non-crossing
half-permutation with zero open blocks. We convert each of
these $k$ half-permutations with one open block to a
half-permutation with zero open blocks by closing the single
open block. Note that this construction is reversible. Since
we always require at least one block we subtract 1 from
$P_0(z)$ and since we convert an open block to a closed
block we divide by $c$.

\end{remark}

\section{Recursion Formula for non-crossing circular
half-permutations}\label{recursion}

In section \ref{firstkind} we constructed a family $\{ \wt
\Gamma_n \}_n$ of polynomials we called shifted Chebyshev
polynomials of the first kind. In the notation established
there $\wt \Gamma_n(x) = \sum_{k = 0}^n g_{n, k}' x^k$ where
$g_{n,k}' \in {\mathbb Z}[c]$. We established the recursion
formula
$$
g_{n,k-1}' = g_{n+1, k}' + (1 + c) g_{n,k}' + c\, g_{n-1,k}'
$$
for $n > 1$ and
$$
g_{1, k-1}' = g_{2,k}' + (1 + c) g_{1, k}' + 2 c g_{0, k}
$$
where $g_{n,k}' = 0$ whenever $n < 0$, $k < 0$ or $k > n$. We
let $\wt \Gamma$ be the matrix with $(n,k)$ entry $g_{n,k}'$ and
$\{g _{n,k} \}_{n,k}$ be the entries of $\wt \Gamma^{-1}$. We
showed that
$$
g_{n+1, k} = g_{n,k-1} + (1 + c) g_{n,k} + c g_{n,k+1} \leqno
(\ref{geq3})
$$
and 
$$
g_{n+1, 0} = (1 + c) g_{n, 0} + 2 c g_{n,1}
\leqno (\ref{geq4})
$$
We wish to show that these relations have an interpretation
in terms of non-crossing circular
half-permutations. Moreover this interpretation is essential
for the proof of our main theorem (Theorem \ref{mt}).

Let $\pi$ be a non-crossing circular half-permutation on $[n]$
and $\#(\pi)_{\cl}$ the number of closed blocks. We wish to
show that for $k \geq 1$ and all $n$ 
\begin{equation}\label{recursioneq1}
g_{n,k} = \sum_\pi c^{ \#(\pi)_{\rm cl} }
\end{equation}
where the sum is over all non-crossing circular
half-permutations on $[n]$ with $k$ open blocks and
$\#(\pi)_{\cl}$ is the number of closed blocks.

Let us denote be $\ncc(n)$ the collection of non-crossing
circular half-permutations on $[n]$ and $\ncc(n)_k$ the
collection of non-crossing circular half-permutations on
$[n]$ with $k$ open blocks.

In order to prove equation (\ref{recursioneq1}), even just
for $k \geq 1$, we have to extend the definition of
non-crossing circular half-permutations (Definition
\ref{circulardef}) to the case $k = 0$.

\begin{definition}\label{circular0}
A non-crossing circular half-permutation on $[n]$ with zero
open blocks is a pair $(\pi, B)$ where $\pi$ is a
non-crossing permutation on $[n]$ and $B$ is either a block
of $\pi$ or a block of $\pi^c$. The block $B$ is called the
{\it designated block}. The collection of non-crossing
circular half-permutations with zero open blocks is denoted
$\ncc(n)_0$.
\end{definition}

\begin{notation}
For all $n \geq 1$ and $n \geq k \geq 1$ let
$$
\ol g_{n,k} = \sum_{\pi \in \ncc(n)_k} c^{\#(\pi)_{\cl} }
$$
and 
$$
\ol g_{n,0} = 
\mathop{\sum_{(\pi, B) \in \ncc(n)_0}}_{B \in \pi^c} c^{\#(\pi)} 
+
\mathop{\sum_{(\pi, B) \in \ncc(n)_0}}_{B \in \pi} c^{\#(\pi) - 1}
$$
For convenience we set $\ol g_{0,0} = 1$. 
\end{notation}

Our goal in this section is to prove

\begin{theorem}\label{gthm}
$\ol g_{n,k} = g_{n,k}$ for $0 \leq k \leq n$
\end{theorem}

We shall prove this theorem by showing that $\{ \ol g_{n,k}
\}_{n,k}$ satisfy the same recurrence relation as $\{g_{n,k}
\}_{n,k}$, namely equations (\ref{geq3}) and (\ref{geq4}).

We will also need to extend the definition of initial and
final points to all the blocks of a half-permutation. Recall
that in Definition \ref{initial} we supposed that $\ol B$
was a block of $\pi^c$ and we defined for an open block $B$
of $\pi$ two points (which coincide if $B$ is a singleton)
which we called the initial and final points. We wish to
extend this to the other blocks of $\pi$ and to the blocks
of a half-permutation with zero open blocks. In the
following definition $\gamma$ is the permutation of $[n]$
with the one cycle $(1, 2, 3, \dots , n)$.

\begin{definition}\label{initialex}
Let $\pi \in NC(n)$ be a non-crossing permutation of $n$ and $B$ 
a block of $\pi$.  
\begin{enumerate}
\item Suppose $\ol B$ a block of $\pi^c$ and $B \cap \ol B =
\emptyset$. Let $j$ be any point of $\ol B$ and let $k$ be
such that $\gamma^k(j) \in B$ but $\gamma^l(j) \not \in B$
for all $1 \leq l < k$. Then $i = \pi^k(j)$ is the {\it
initial point} of $B$ (relative to $\ol B$) and
$\pi^{-1}(i)$ is the {\it final point} of $B$ (relative to
$\ol B$). By the non-crossing property of $\pi$, $i$ is
independent of the choice of $j$.

\item Suppose $B_1 \not = B$ is a block of $\pi$. We choose
any $j \in B$ and find the first $k$ such that $\gamma^k(j)
\in B_1$, i.e. $\gamma^k(j) \in B_1$ but $\gamma^l(j) \not
\in B_1$ for $1 \leq l < k$. Let $i = \gamma^k(j)$, we call
$i$ the {\it initial point} of $B_1$ (relative to $B$) and
$\pi^{-1}(i)$ the {\it final point} of $B_1$. Again $i$ is
not affected by the choice of $j$.
\end{enumerate}
\end{definition}

\begin{remark}
If $\ol B$ is a block of $\pi^c$ and $B$ is a block of $\pi$
with $B \cap \ol B \not = \emptyset$, then Definition
\ref{initial} and ({\it a}) above coincide.
\end{remark}

In order to prove Theorem \ref{gthm} we introduce a new way of describing
elements of $\ncc(n)_k$ using a dot structure. By a dot structure we
mean placing either black or white dots on each of $1, 1', 2, 2',
\dots , n, n'$. 

\begin{notation}
Let $D_{j,k,n}$ be the collection of dot structures on $1,
1', 2, 2',\ab \dots , n, n'$ such that
\begin{itemize}
\item there are $j$ black dots on primed numbers and the remaining $(n - j)$
primed numbers have white dots;
\item there are $(k + j)$ white dots on unprimed numbers and black dots on the
remaining $n - (j + k)$ unprimed numbers.
\end{itemize}

\end{notation}

\begin{theorem}\label{dotthm1}
For $k > 0$ there is a bijection between the subset $\{ \pi
\in \ncc(n)_k \mid \pi$ has $j$ closed blocks $\}$ and
$D_{j,k,n}$.
\end{theorem} 

\begin{proof}
Suppose we are given a dot structure.  We place the numbers
$1, 1', 2, 2', \dots , n, n'$ around the circle in a
clockwise fashion.

\bigskip\noindent
$\vcenter{\hsize140pt\includegraphics{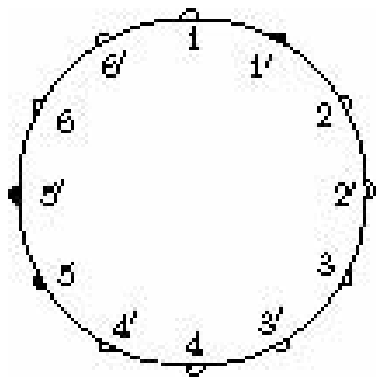}}$
$\vcenter{\hsize190pt\parindent0pt\small\raggedright {\bf
Figure \figno.} In this example $n= 6$, $j = 2$, and $k =
3$. We place black dots on $1'$ and $5'$ and white dots on
$2'$, $3'$, $4'$ and $6'$. $j + k = 5$ so we place 5 white
dots on a unprimed numbers, 1, 2, 3, 4, and 6 say. To
complete we place a black dot on $5'$.}$

Starting at any black dot we move in a counter-clockwise
direction until we come to the first {\it available} white
dot and connect these two dots. By available we mean that
every time we pass over a black dot we must skip over an
additional white dot. Since there are in total $n - k$ black
dots and $n + k$ white dots, every black dot can be
connected to a white dot in the manner described. Moreover
since a black dot on a primed number can only be paired with
a white dot on an unprimed number, there will be $k$
remaining white dots on unprimed numbers. Similarly a black
dot on an unprimed number can only be paired with a white
dot on a primed number so there will be $k$ white dots on
primed numbers remaining.

\noindent
$\vcenter{\hsize160pt\includegraphics{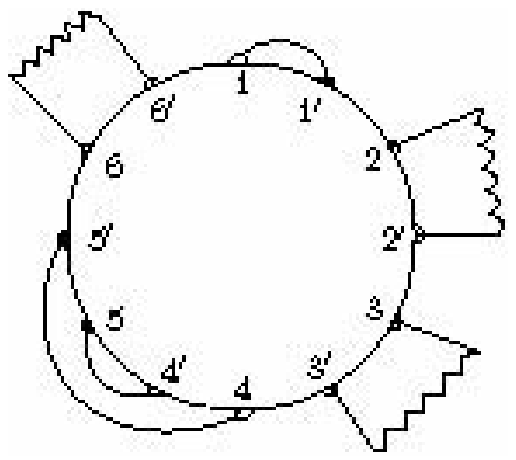}}$\quad
$\vcenter{\hsize160pt\parindent0pt\small\raggedright {\bf
Figure \figno.} We then connect each black dot to its first
available neighbour in the counter-clockwise direction. The
$k$ remaining white dots on primed numbers get connected to
the $k$ remaining white dots on unprimed numbers.}$
\medskip

We connect each remaining white dot on a primed number to
the first available white dot on an unprimed number, moving
as before in the counter-clockwise direction.

\smallskip\hfill
$\vcenter{\hsize160pt\includegraphics{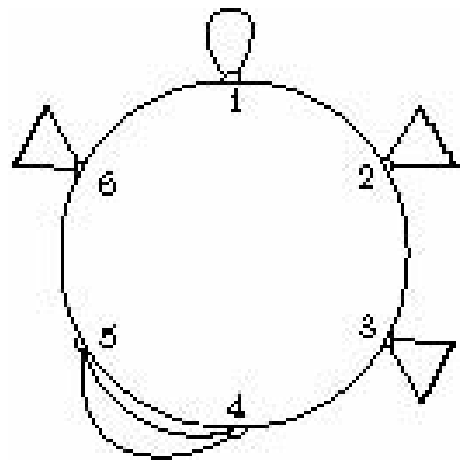}}
\vcenter{\hsize160pt\parindent0pt\small\raggedright
{\bf Figure \figno.} We obtain a non-crossing circular
half-permutation  on $[6]$ with 2 closed blocks and 3 open
blocks.}$\hfill\hbox{}

\smallskip%\hrule\medskip

Finally we squeeze together each unprimed number $i$ with
its corresponding primed number $i'$ to form a non-crossing
circular half-permutation. Each black dot on a primed number
starts a closed block. There are $n - j$ white dots on
primed numbers of which $n - j - k$ are paired with black
dots on unprimed numbers.

Each white dot on a primed number starts an open block and
so we obtain $k$ open blocks and $j$ closed blocks. This
incidentally shows that the number of such permutations is
${ n \choose j}{ n \choose j + k}$; so the coefficient of
$c^j$ in
$$
\ds\sum_{\pi \in \ncc(n)_k} c^{
\#(\pi)_{\rm cl} } \mbox{ is  } { n \choose j}{ n \choose j + k}.
$$ 
 
This construction is reversible: starting with a
non-crossing circular half-permutation with $k$ open blocks
and $j$ closed blocks we arrange black and white dots as
follows.

\medskip\noindent $\vcenter{\hsize160pt
\includegraphics{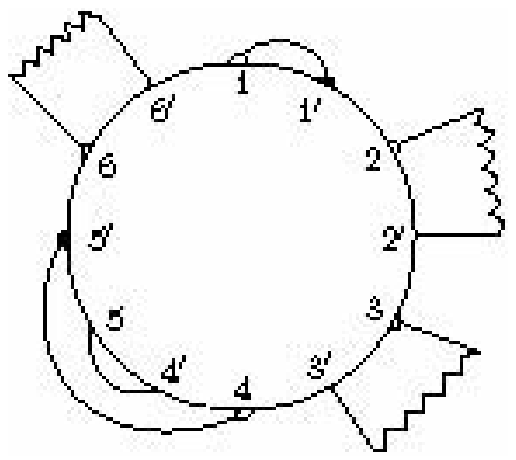}}$\qquad
$\vcenter{\hsize150pt\small\raggedright\small {\bf Figure
\figno.} We start with a non-crossing circular
half-permutation and we place a black dot on $1'$ and
$5'$. We place white dots on 1, 2, 3, 4, and 6.}$

On the unprimed numbers we place a white dot on the initial
point of each block and black dots on the $n - j - k$
remaining points. On the primed numbers we place a black dot
on the final point of each of the $j$ closed block and a
white dot on each of the remaining points. This establishes
a bijection between $\ncc(n)_k$ and $D_{j,k,n}$.
\end{proof}
\noindent
$\vcenter{\hsize160pt \includegraphics{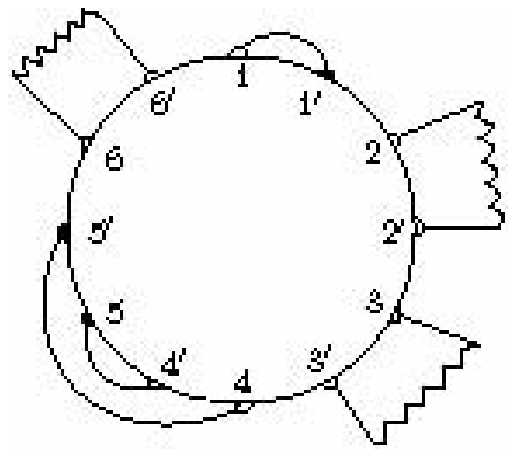}}$\qquad
$\vcenter{\hsize150pt\small\raggedright\small {\bf Figure
\figno.} We then place a black dot on all the remaining
unprimed numbers and a white dot on all the remaining primed
numbers.}$

\begin{theorem}
There is a bijection between $D_{j,0,n}$ and the subset of
$\ncc(n)_0$: $\{ (\pi, B) \mid \pi$ has $j$ blocks when $B
\in \pi^c$, or $\pi$ has $j+1$ blocks when $B \in \pi \}$.
\end{theorem}

\begin{proof}
Let $(\pi, B) \in \ncc(n)_0$ with $B \in \pi^c$. Suppose
that $\pi$ has $j$ blocks, we shall construct a dot diagram
in $D_{j, 0, n}$. We put a black dot on the prime of the
final point of each block of $\pi$, all other primed numbers
get a white dot. We put a white dot on the unprime of the
initial point of each block of $\pi$; all other unprimes get
a black dot. If a block of $\pi$ is a singleton then we put
a black dot on the prime and a white dot on the unprime of
the number. Thus there are $j$ black dots on primed numbers
and $j$ white dots on unprimed numbers; so the dot diagram
is in $D_{j,0,n}$.

Next suppose that $(\pi, B) \in \ncc(n)_0$ with $j+1$ blocks
and $B \in \pi$. We follow the same construction as before
for each of the $j$ undesignated blocks. For the designated
block $B$ we put a white dot on each unprimed number and a
black dot on each primed number.  Thus on the designated
block there are no black dots on primed numbers or white
dots on unprimed numbers. So again we obtain a dot pattern
in $D_{j,0,n}$.

Now suppose we have a dot pattern in $\ncc(n)_0$ we must
construct a non-crossing permutation $\pi$ and identify
block $B$ either of $\pi$ or of $\pi^c$.

We obtain $\pi$ as in the proof of Theorem \ref{dotthm1}. We
have $n$ black dots and $n$ white dots. We connect each
black dot to the first available white dot. The blocks of
$\pi$ are constructed as follows. Start at any of the $j$
white dots on an unprimed number and `open' a block. We move
in the clockwise direction adding points until we come to a
black dot on a primed number, where we `close' the block. If
we encounter another white dot on an unprimed number before
closing the block we open new block and close this one
before closing the first one opened. These Catalan paths can
only produce non-crossing permutations.

If this matching exhausts all of the points then we are in
the case $B \in \pi^c$. We can identify the points of $B$ as
follows. A point $i$ is in $B$ if it passes the following
test. Starting at the unprime of $i$ and moving in a
clockwise direction we count the number of white points
passed minus the number of black points. If we can return to
$i$ without the running total ever being negative, then $i
\in B$.

If the matching doesn't exhaust all of the points then there
are points remaining of the form black dot on unprimed
number and white dot on unprimed number. These will then
form the points of the designated block $B \in \pi$. As we
can only encounter a black-white pair when there are no open
blocks, $B$ cannot cross any of the previously constructed
blocks. Thus $(\pi, B) \in \ncc(n)_0$.\end{proof}

\begin{theorem}
For $k > 0$
\begin{equation}\label{gbar1}
\ol g_{n+1,k} = \ol g_{n,k-1} + (1 + c) \ol g_{n,k} + c \ol g_{n,k+1}
\end{equation}
and for $k = 0$
\begin{equation}\label{gbar2}
\ol g_{n+1,0} = (1 + c) \ol g_{n,0} + 2 c \ol g_{n, 1}
\end{equation}
\end{theorem}

\begin{proof}
For all $n$ and $k$ we write $\ncc(n)_k$ as the disjoint
union of four subsets according to the four possible dot
patterns on $n$:
\begin{eqnarray*}
\ncc(n)_{k,1} &=& \{ (\pi, B) \in \ncc(n)_k \mid 
\mbox{ white on } n \mbox{ white on } n' \}\\
\ncc(n)_{k,2} &=& \{ (\pi, B) \in \ncc(n)_k \mid 
\mbox{ black on } n \mbox{ white on } n' \}\\
\ncc(n)_{k,3} &=& \{ (\pi, B) \in \ncc(n)_k \mid 
\mbox{ white on } n \mbox{ black on } n' \}\\
\ncc(n)_{k,4} &=& \{ (\pi, B) \in \ncc(n)_k \mid 
\mbox{ black on } n \mbox{ black on } n' \}\\
\end{eqnarray*}

We shall define maps
\begin{eqnarray*}
&&\phi_1 : \ncc(n+1)_{k,1} \rightarrow \ncc(n)_{k-1} \\
&&\phi_2 : \ncc(n+1)_{k,2} \rightarrow \ncc(n)_{k} \\
&&\phi_3 : \ncc(n+1)_{k,3} \rightarrow \ncc(n)_{k} \\
&&\phi_4 : \ncc(n+1)_{k,4} \rightarrow \ncc(n)_{k+1} \\
\end{eqnarray*}
such that $\#(\phi_i(\pi))_{\cl} = \#(\pi)_{\cl}$ for $i = 1, 2$ and
$\#(\phi_i(\pi))_{\cl} = \#(\pi)_{\cl} - 1$ for $i = 3, 4$. 

On the level of dot diagrams the maps $\phi_i$ all simply
delete the points $n + 1$, $(n + 1)'$ and their dots except
in the case $k = 0$ and $i = 1$. When $k = 0$ and $i = 1$,
$\phi_1$ deletes the two white dots on $n+1$ and $(n+1)'$
and then reverses the colour of all dots.

Let us check that each of the maps $\phi_i$ is a bijection
and has the required effect on the number of closed blocks.

{\it The case} $i = 1$. We remove a white dot on an unprimed
number so the number of blocks decreases by one, however the
number of black dots on primed numbers is unchanged so the
block lost is a white block. Thus $\#(\phi_1(\pi))_{\cl} =
\#(\pi)_{\cl}$. Since we had $k$ open blocks to start with
$\phi_1(\pi) \in \ncc(n)_{k-1}$. The inverse map from
$\ncc(n)_{k-1}$ is to insert $n+1$ and $(n+1)'$ and put
white dots on them.

{\it The case} $i = 2$. We remove a black dot on an unprimed
number so the number of blocks is unchanged, moreover the
number of black dots on primed numbers is unchanged so the
number of closed blocks is unchanged. Thus
$\#(\phi_2(\pi))_{\cl} = \#(\pi)_{\cl}$. Since we had $k$
open blocks to start with $\phi_2(\pi) \in \ncc(n)_{k}$. The
inverse map from $\ncc(n)_{k-1}$ is to insert $n+1$ and
$(n+1)'$ and put a black dot on $n + 1$ and a white dot on
$(n+1)'$.

{\it The case} $i = 3$. We remove a white dot on an unprimed
number so the number of blocks is decreased by one, also the
number of black dots on primed numbers is decreased by one
so the number of closed blocks is decreased by one. Indeed,
in this case we see that $n +1$ was a closed singleton and
we removed it. Thus $\#(\phi_3(\pi))_{\cl} = \#(\pi)_{\cl} -
1$. Since we had $k$ open blocks to start with $\phi_3(\pi)
\in \ncc(n)_{k}$. The inverse map from $\ncc(n)_{k-1}$ is to
insert $n+1$ and $(n+1)'$ and put a white dot on $n + 1$ and
a black dot on $(n+1)'$.

{\it The case} $i = 4$. We remove a black dot on an unprimed
number so the number of blocks is unchanged, however the
number of black dots on primed numbers is reduced by one so
a closed block is lost and an open block is gained.  Thus
$\#(\phi_4(\pi))_{\cl} = \#(\pi)_{\cl} - 1$. Since we had
$k$ open blocks to start with $\phi_4(\pi) \in
\ncc(n)_{k+1}$. The inverse map from $\ncc(n)_{k+1}$ is to
insert $n+1$ and $(n+1)'$ and put black dots on them.

This construction works for all $n$, $k$, and $i$ except
when $k = 0$ and $i = 1$. In this situation the number of
black dots equals the number of white dots (because $k = 0$)
and we remove two white dots leaving two fewer white dots
than black, which is forbidden. So in this case we remove
the two white dots and change all black dots to white and
all white to black. So when $k = 0$ both $\phi_1$ and
$\phi_4$ map into $\ncc(n)_{1}$. This justifies the 2 in
equation (\ref{gbar2}).
\end{proof}

\begin{theorem}\label{circularrecursion}
For $n \geq 1$ and $k \geq 1$
$$
g_{n,k} = \sum_{\pi \in \ncc(n)_k} c^{\#(\pi)_{\cl} }
$$
and
$$
g_{n,0} = \sum_{\pi \in NC(n)} \#(\pi^c) c^{\#(\pi) } + \#(\pi) c^{ \#(\pi) - 1}
$$
\end{theorem}

\begin{proof}
We only have to check that we have equality of $\{ g_{1,0},
g_{1,1} \}$ and $\{ \ol g_{1,0}, \ol g_{1,1} \}$ so that the
induction can start. We have by calculation $g_{1,0} = 1 +
c$ and $g_{1,1} = 1$. On the level of diagrams when $n = 1$
we have one block with one element in it. When $k = 0$ we
can either choose $B$ to be this block or the single block
in the complement so that $\ol g_{1,0} = 1 + c$. When $k =
1$ we must have one open block and since there is only one
block, zero closed blocks, thus $\ol g_{1,1} = 1$.
\end{proof}

\begin{remark}\label{circulargenfunction}
There is an intriguing interpretation of equation
(\ref{power1}) for non-crossing circular half-permutations
analogous to that (remark \ref{lineargenremark}) for
non-crossing linear half-permutations.
$$ G_n(z) = \bigg(\frac{P_0(z)
-1}{c}\bigg)^n G_0(z) \leqno (9)$$
Each non-crossing circular half-permutation with $k$ open blocks
can be assembled from $k$ non-crossing linear half-permutations
with one open block and one non-crossing circular
half-permutation with zero open blocks. The point of insertion
of linear permutations is determined by the \textit{designated}
block. We shall not present the details here.
\end{remark}

\section{The Proof of Theorem \ref{mt}}

Let for each $N$, $X_{N,1}, \dots X_{N,p}$ be $p$
independent $N \times N$ Wishart random matrices
i.e. $X_{N,i} = G_{N,i}^\ast G_{N,i}$ with $G_{N,i} = \{
g^{(N,i)}_{m,n} \}_{m,n}$ an $M \times N$ Gaussian random
matrix with independent entries and $\E( |g^{(n,
i)}_{m,n}|^2 ) = 1/N$. Fix $k,l > 1$ and let $\vec m = (m_1,
\dots , m_k)$ and $\vec n = (n_1, \dots , n_l)$ with $m_s >
0$ and $n_t > 0$ for all $s$ and $t$.  Suppose also that we
have $i_1, \dots i_k$ and $j_1, \dots j_l$ with $1 \leq i_s,
j_t \leq p$ which are each {\it cyclically alternating}
i.e. for $|s_1 - s_2| = 1$, $i_{s_1} \not = i_{s_2}$ and for
$| t_1 - t_2| = 1$, $j_{t_1} \not = j_{t_2}$ and in addition
$i_k \not = i_1$ and $j_l \not = j_1$.

Let
$$
S_{\vec m, \vec i} = \Tr(\Pi_{m_1}(X_{N,i_1})\ab \cdots
\Pi_{m_k}(X_{N,i_k})) $$
and 
$$
S_{\vec n, \vec j} = \Tr(\Pi_{n_1}(X_{N, j_1}) \cdots
\Pi_{n_l}(X_{N,j_l})) $$

From \cite{mn} Theorem 7.5
we have.
\begin{theorem}\label{mnthm7.5}
\[
\lim_N \kappa_2( \Tr( X_{N,i_1}^{m_1} \cdots X_{N,i_k}^{m_k} ),
                 \Tr( X_{N,j_1}^{n_1} \cdots X_{N,j_l}^{n_l} ) )
= |S_{NC}( \vec m, \vec n)|_c
\]
where the undefined notation is explained below.
\end{theorem}

\begin{notation}
Let $\vec u = (u_1, \dots , u_k)$ and $\vec v = (v_1, \dots
, v_l)$ where $u_1, u_2,\ab \dots , u_k$ and $v_1, v_2,
\dots , v_l$ are positive integers. Let $u = u_1 + u_2 +
\cdots + u_k$ and $v = v_1 + v_2 + \cdots + v_l$.

Let $S_{NC}(\vec u, \vec v)$ be the set
of non-crossing annular $(u,v)$-permutations where 
\begin{itemize} 
\item on the $u$-circle we have $k$ intervals $I_1$, $I_2$,
      \dots , $I_k$ of points with $u_s$ points in the interval
      $I_s$ all of colour $i_s$; 
\item on the $v$-circle we have $l$ intervals $J_1$, $J_2$,
      \dots , $J_l$ of points with $j_t$ points in the interval
      $J_t$ all of colour $j_t$;
\item only points of the same colour can be connected.

\end{itemize}
To simplify the notation we have omitted the dependence on
the colours $i_1, \dots , i_k$ and $j_1, \dots , j_l$; since
we only work with one set of colours this should not cause
any confusion.
\end{notation}

\begin{notation}
For $\vec x = (x_1, \dots , x_k)$ and $\vec y = (y_1, \dots
, y_l)$ let $S_{NC}\big( \myatop{\vec u}{\vec x},
\myatop{\vec v }{ \vec y}\big)\ab = S_{NC}(\myatop{u_1 }{
x_1}, \dots, \myatop{u_k }{ x_k};\ab \myatop{v_1 }{ y_1},
\dots , \myatop{v_l }{ y_l}\big) $ be the permutations in
$S_{NC}( \vec u, \vec v)$ such that from the interval $I_s$
there are exactly $x_s$ blocks that meet the $v$-circle {\it
and} from the interval of $J_t$ there are exactly $y_t$
blocks that meet the $u$-circle, i.e. $I_s$ has $x_s$
through-blocks and $J_t$ has $y_t$ through-blocks.

If $x_s = u_s$ then we shall write
$S_{NC}(\myatop{u_1}{x_1}, \dots, \underline{u}_s , \dots,
\myatop{u_k }{ x_k};\ab \myatop{v_1 }{ y_1}, \dots ,
\myatop{v_l }{ y_l}\big) $ in place of $S_{NC}(\myatop{u_1
}{ x_1}, \dots, \myatop{u_s }{ u_s}, \dots , \myatop{u_k }{
x_k};\ab \myatop{v_1 }{ y_1}, \dots , \myatop{v_l }{
y_l}\big) $. Note that the condition $x_s = u_s$ means that
every point of $I_s$ is in a different block and each of
these blocks is a through-block.

If $\vec x = \vec u$ and $\vec y = \vec v$ we shall write
$S_{NC}(\vec{\underline{u}}, \vec{\underline{v}})$ in place
of $S_{NC}\big( \myatop{ \vec u }{ \vec u}, \myatop{\vec v
}{ \vec v}\big)$. In this case every block is a pair
connecting the two circles (i.e. a spoke diagram),
respecting the colours.

\end{notation}

\begin{notation}
Suppose $A \subset S_{NC}(m,n)$. Let $|A|_c = \sum_{\pi \in
A} c^{\#(\pi)}$. Suppose $1 \leq s < t \leq k$. Let
$NC([s,t])$ be the non-crossing partitions $\pi$ of $[s, t]
= \{s, s+1, \dots , t\}$ that respect the colouring given by
$i_s, \dots , i_t$, i.e. are such that if $a$ and $b$ are in
the same block of $\pi$ then $i_a = i_b$.

Let $X_{s,t} = \cup_{s \leq r \leq t} I_r$ with the
colourings given by $i_s, \dots , i_t$, i.e. the points of
$I_r$ are coloured $i_r$. Let $NC(X_{s,t})$ be the
non-crossing partitions of $X_{s,t}$ that respect the
colouring, i.e. each block only contains one colour.

Given $\pi \in NC(X_{s,t})$ we obtain a partition $\nu_\pi
\in NC([s,t])$ in which $r_1, r_2 \in [s,t]$ are in the same
block whenever a block of $\pi$ meets both the intervals
$I_{r_1}$ and $I_{r_2}$. The idea is to take a partition of
$X_{s,t}$ that respects the colours and only remember how
the partition connects the intervals, i.e. shrink each
interval to a point.

For $B \subset [s, t]$ such that $i_{r_1} = i_{r_2}$ for all
$r_1, r_2 \in B$, let $X_B = \cup_{r \in B} I_r$. Let
$\sigma_B$ be the partition of $X_B$ in which each interval
$I_r$ $(r \in B)$ is a block of $\sigma_B$.

Let 
$$\ds p_B = 
\mathop{\sum_{\pi \in NC(X_B)}}_{\sigma_B \vee \pi = 1_{X_B}}
c^{\#(\pi)}
$$ i.e the sum is over all $\pi$ that connect all the
intervals $\{I_r \}_{r \in B}$.
\end{notation}

\begin{remark}\label{remark16}
{\it (i)} Suppose $s < t$, $x_s = \cdots = x_t =0$, $x_{s-1}
> 0$, $x_{t+1} > 0$ and $i_{s-1} \not= i_{t+1}$.  Then by
Theorem \ref{decomposition},
$$ S_{NC}\Big(\myatop{\vec u }{ \vec x} ; \myatop{\vec v }{
\vec y}\Big) = S_{NC}\Big(\myatop{u_1 }{ x_1} , \dots ,
\myatop{u_{s-1} }{ x_{s-1} }, \myatop{u_{t+1} }{ x_{t+1} },
\dots , \myatop{u_k }{ x_k}; \myatop{\vec v }{ \vec y}\Big)
\times NC(X_{s,t})
$$ If $i_{s-1} = i_{t+1}$ then
$$S_{NC}\Big(\myatop{\vec u }{ \vec x} ; \myatop{\vec v }{ \vec y}\Big)
= S_{NC}\Big(\myatop{u_1 }{ x_1} , \dots , \myatop{u_{s-1} + u_{t+1} }{
x_{s-1} + x _{t+1}}, \dots , \myatop{u_k }{ x_k};
\myatop{\vec v }{ \vec y}\Big) \times NC(X_{s,t})$$

{\it (ii)} If $B = \{w \}$ is a singleton then 
$$ p_B = \sum_{\pi \in NC(I_w)} c^{\#(\pi)} = p_{u_w,0}
$$
\end{remark}

\begin{lemma}\label{lemma17}
$$ | NC( X_{s,t} )|_c = \mathop{\sum_{\tau \in NC([s,t]
    )}}_{\tau = \{B_1, \dots , B_q\}} \prod_{1 \leq i \leq
  q} p_{B_i}
$$

\end{lemma}%\hfill\framebox[7pt]{\relax}

\begin{proof}
We partition $NC(X_{s,t})$ over $NC([s, t])$, i.e. write
$$ NC(X_{s,t}) = \ds\bigcup_{\tau \in NC([s,t])} \{ \pi \in
NC(X_{s,t}) \mid \nu_\pi = \tau \}
$$ where $\nu_\pi$ is the partition of the colours induced
by $\pi$. Combine this with the observation that for $\tau =
\{B_1, \dots , B_q \}$
$$\displaylines{
\{ \pi \in NC(X_{s,t}) \mid \nu_\pi = \tau \}
=
\{ \rho \in NC(X_{B_1}) \mid \rho \vee \sigma_{B_1} = 1_{X_{B_1}}
\} \hfill\cr
\times \cdots \times 
\{ \rho \in NC(X_{B_q}) \mid \rho \vee \sigma_{B_q } =
1_{X_{B_q}} \} 
\cr}$$
\end{proof}

\begin{lemma} \label{lemma18}
 Suppose that $x_i = 0$ for some  $i$. Then 
$$ 
\mathop{\sum_{u_1, \dots , u_k}}_{v_1, \dots , v_l}
\mathop{\prod_{1 \leq r \leq k}}_{1 \leq s \leq l} p_{m_r,
u_r}' p_{n_s, v_s}' \Big|S_{NC}\Big( \myatop{\vec m }{ \vec
x}, \myatop{\vec n }{ \vec y}\Big)\Big|_c = 0
$$
\end{lemma}

\begin{proof} 
By invariance of the trace under cyclic permutatione we may
suppose $x_1 = \cdots = x_t = 0$ and $x_{t+1} > 0$.
Furthermore we shall assume that $i_{t+1} \not= i_k$, the
proof when $i_{t+1} = i_k$ is handled similarly.

$$\displaylines{ \mathop{\sum_{u_1, \dots , u_k}}_{v_1,
\dots , v_l} \mathop{\prod_{1 \leq r \leq k}}_{1 \leq s \leq
l} p_{m_r, u_r}' p_{n_s, v_s}' \Big|S_{NC}\Big( \myatop{ m_1
}{ 0}, \dots , \myatop{ m_{t} }{ 0}, \myatop{m_{t+1} }{
x_{t+1}} , \dots , \myatop{m_k }{ x_k}, \myatop{\vec n }{
\vec y}\Big)\Big|_c \hfill \cr = \mathop{\sum_{u_1, \dots ,
u_k}}_{v_1, \dots , v_l} \mathop{\prod_{1 \leq r \leq k}}_{1
\leq s \leq l} p_{m_r, u_r}' p_{n_s, v_s}'
\Big|S_{NC}\Big(\myatop{m_{t+1} }{ x_{t+1}} , \dots ,
\myatop{m_k }{ x_k}, \myatop{\vec n }{ \vec y}\Big)\Big|_c
\times | NC( X_{0,t})|_c \hfill \cr}$$

By Lemma \ref{lemma17} we can write $|NC(X_{0,t})|_c$ as a
sum over $NC([0,t])$. Since the colours are cyclically
alternating we know that each $\tau$ in $NC([0,t])$ must
contain at least one singleton $\{a(\tau)\}$ and by Remark
\ref{remark16} ({\it ii}\/) the corresponding $p_B$ is a
$p_{u_{a(\tau)}, 0}$; we have that for some $M_\tau$
$$
| NC(X_{0,t}) |_c =
\sum_{\tau \in NC([0,t])}
p_{u_{a(\tau)},0}\, M_\tau
$$

Inserting this into our equation above we continue

\begin{align*}
&= \sum_{\tau \in NC([0,t])} \mathop{\sum_{u_i, \dots ,
u_k}}_{v_1, \dots , v_l} \mathop{\prod_{1 \leq r \leq k}}_{1
\leq s \leq l} p_{m_r, u_r}' p_{n_s, v_s}' p_{u_{a(\tau)},0}
M_\tau \Big|S_{NC}\Big(\myatop{m_{t+1} }{ x_{t+1}} , \dots ,
\myatop{m_k }{ x_k}, \myatop{\vec n }{ \vec y}\Big)\Big|_c
\hfill \\ &= 0
\end{align*}

Because, by Theorem \ref{linearrecursion},
$\sum_{u_{a(\tau)}} p_{m_{a(\tau)}, u_{a(\tau)}}'
p_{u_{a(\tau)},0} = \delta_{m_{a(\tau),0}} = 0$.
\end{proof}

\begin{proposition}\label{bigprop}
$$ \lim_{N \rightarrow \infty} \kappa_2(S_{\vec m, \vec i},
S_{\vec n, \vec j}) = | S_{NC}(\vec{\underline m},
\vec{\underline n})|_c
$$
\end{proposition}

\begin{proof}
$$\displaylines{ \lim_N \kappa_2(S_{\vec m, \vec i}, S_{\vec
n, \vec j}) \hfill \cr = \lim_N \mathop{\sum_{u_1, \dots ,
u_k}}_{v_1, \dots , v_l} \mathop{\prod_{1 \leq r \leq k}}_{1
\leq s \leq l} p_{m_r, u_r}' p_{n_s, v_s}' \kappa_2(
X_{N,i_1}^{u_1} \cdots X_{N,i_k}^{u_k}, X_{N,j_1}^{v_1}
\cdots X_{N,j_l}^{v_l} ) \hfill\cr
\stackrel{Thm. \ref{mnthm7.5}}{=} \mathop{\sum_{u_1, \dots ,
u_k}}_{v_1, \dots , v_l} \mathop{\prod_{1 \leq r \leq k}}_{1
\leq s \leq l} p_{m_r, u_r}' p_{n_s, v_s}' \Big| S_{NC}(
u_1, \dots , u_k; v_1, \dots , v_l )\Big|_c \hfill \cr}$$

    {\narrower we write $S_{NC}( u_1, \dots , u_k; v_1, \dots ,
    v_l )$ as a disjoint union over all $\vec x$ and $\vec y$ of
    $S_{NC} \Big( \myatop{u_1 }{ x_1},\ab \dots , \myatop{u_k }{
    x_k}; \myatop{v_1 }{ y_1} , \dots , \myatop{v_l }{ y_l}
    \Big)$ \par}

$$\displaylines{ = \mathop{\sum_{u_1, \dots , u_k}}_{v_1,
\dots , v_l} \mathop{\prod_{1 \leq r \leq k}}_{1 \leq s \leq
l} p_{m_r, u_r}' p_{n_s, v_s}' 
\sum_{\vec x, \vec y} \Big| S_{NC}
\Big( \myatop{u_1 }{ x_1}, \dots , \myatop{u_k }{ x_k};
\myatop{v_1 }{ y_1} , \dots , \myatop{v_l }{ y_l}
\Big)\Big|_c \hfill\cr}$$ 

{\narrower by Lemma \ref{lemma18} we have,\par}

$$\displaylines{ = \mathop{\sum_{u_1, \dots , u_k}}_{v_1,
\dots , v_l} \mathop{\prod_{1 \leq r \leq k}}_{1 \leq s \leq
l} p_{m_r, u_r}' p_{n_s, v_s}' \mathop{\sum_{\vec x, \vec
y}}_{x_1, \dots , x_k > 0 } 
\Big| S_{NC} \Big(
\myatop{u_1 }{ x_1}, \dots , \myatop{u_k }{ x_k};
\myatop{v_1 }{ y_1} , \dots , \myatop{v_l }{ y_l}
\Big)\Big|_c \hfill\cr}$$ {\narrower by Lemma \ref{lemma18}
again we have, \par}
$$\displaylines{
= 
\mathop{\sum_{u_1, \dots , u_k}}_{v_1, \dots , v_l}
\mathop{\prod_{1 \leq r \leq k}}_{1 \leq s \leq l}
p_{m_r, u_r}' p_{n_s, v_s}'  
\mathop{ 
\sum_{x_1, \dots , x_k > 0 } }_{y_1, \dots , y_l > 0 } 
\Big| S_{NC} \Big( \myatop{u_1 }{ x_1}, \dots , \myatop{u_k }{ x_k}; 
\myatop{v_1 }{ y_1} , \dots , \myatop{v_l }{ y_l} \Big)\Big|_c \hfill
\cr }$$
{\narrower
Now we have for each interval at least one block passing to the
opposite circle. Thus on each interval we have a non-crossing
linear half-permutation. Thus 
\begin{align*}
 \Big| S_{NC} (\myatop{\vec u }{ \vec x}; \myatop{\vec v }{ \vec
y}\Big) \Big|_c & =
\sum_{\pi \in \ncl(I_1)_{x_1}}
c^{\#(\pi)_{\rm cl}} \ 
\Big| S_{NC} \Big(\myatop{x_1 }{ x_1}, \myatop{u_2 }{ x_2} , \dots ,
\myatop{u_k }{ x_k} ; \myatop{\vec v }{ \vec y} \Big) \Big|_c \\
&=
p_{u_1,x_1} \Big| S_{NC} \Big(\underline{x_1}, \myatop{u_2 }{ x_2}
, \dots , \myatop{u_k }{ x_k} ; \myatop{\vec v }{ \vec y} \Big) \Big|_c 
\end{align*}
by Theorem \ref{linearrecursion}. We can repeat this for each $x_i$
and $y_j$ and thus obtain that our limit
 \par} 
$$\displaylines{ =
\mathop{ 
\sum_{x_1, \dots , x_k > 0 } }_{y_1, \dots , y_l > 0 }
\mathop{\sum_{u_1, \dots , u_k}}_{v_1, \dots , v_l}
\mathop{\prod_{1 \leq r \leq k}}_{1 \leq s \leq l}
p_{m_r, u_r}' p_{n_s, v_s}' p_{u_r, x_r} p_{v_s, y_s} 
| S_{NC} ( \vec{\underline x}; \vec{\underline y} )|_c \hfill \cr 
=
\mathop{ 
\sum_{x_1, \dots , x_k > 0 } }_{y_1, \dots , y_l > 0 }
| S_{NC} ( \vec{\underline x}; \vec{\underline y} )|_c
\mathop{\prod_{1 \leq r \leq k}}_{1 \leq s \leq l}
\delta_{m_r, x_r} \delta_{n_s, y_s} \hfill \cr
= | S_{NC} ( \vec{\underline m}; \vec{\underline n})|_c \hfill
\cr}$$
\end{proof}

\begin{proposition}\label{gsindep}
Let $k > 1$, $1 \leq i \leq p$ and  $i_1, i_2, \dots , i_k$ be cyclically
alternating, and $n, \ab m_1, \dots , m_k > 0$. Then 
$$
\lim_N \kappa_2(S_{\vec m, \vec i},\Gamma_n(X_{N,i})) =0
$$
\end{proposition}

\begin{proof}
As in the proof of Proposition \ref{bigprop}'

\begin{align*}
\lim_N &\kappa_2(S_N,\Gamma_n(X_{N,i})) 
 \stackrel{Thm. \ref{mnthm7.5}}{=} 
\sum_{u_1, \dots , u_k,v}\ 
\prod_{1 \leq r \leq k} p_{m_r, u_r}' q_{n,v}' \ 
\Big| S_{NC} \Big( u_1, \dots ,u_k; v \Big) \Big|_c \\
&=
\sum_{x_1, \dots , x_k > 0}
\sum_{u_1, \dots , u_k,v}\ 
\prod_{1 \leq r \leq k} p_{m_r, u_r}' q_{n,v}' \ 
\Big| S_{NC} \Big( 
\myatop{u_1 }{ x_1}, \dots , \myatop{u_k }{ x_k}; v \Big) \Big|_c \\
&=
\sum_{x_1, \dots , x_k > 0}
\sum_{u_1, \dots , u_k,v}\ 
\prod_{1 \leq r \leq k} p_{m_r, u_r}' p_{u_r, x_r} q_{n,v}' \ 
| S_{NC} ( \vec{\underline x}; v)|_c \cr
&=
\sum_{x_1, \dots , x_k > 0, v}
\prod_{1 \leq r \leq k} \delta_{m_r, x_r} q_{n,v}'
| S_{NC} ( \vec{\underline x}; v)|_c \\
&=
\sum_{1 \leq v \leq n}  q_{n,v}' 
| S_{NC} ( \vec{\underline m}; v)|_c \\ &= 0
\end{align*}

because $S_{NC} ( \vec{\underline m}; v)$ is empty.
\end{proof}

\begin{proposition}\label{ggindep}
Suppose $m, n > 0$. 
$$
\lim_N \kappa_2( \Gamma_m(X_{N,i}), \Gamma_n( X_{N, j}) ) 
= \delta_{i,j} \delta_{m,n}\,  m\, c^m
$$
\end{proposition}

\begin{proof}
$$
\lim_N \kappa_2( \Gamma_m(X_{N, i}), \Gamma_n( X_{N,j} )) 
=
\mathop{\sum_{1 \leq u \leq m}}_{1 \leq v \leq n}
q_{m,u}' q_{n,v}' \lim_N \kappa_2( X_{N,i}^u, X_{N,j}^v )
$$

Now $ \lim_N \kappa_2( X_{N,i}^u, X_{N,j}^v ) = 0$ unless $i = j$, so 
for the rest of the proof we shall assume that $i = j$. Moreover
as $i = j$, $\lim_N \kappa_2( X_{N}^u, X_{N}^v) = | S_{NC}( u; v)
|_c$. Thus
\begin{align*}
\lim_N \kappa_2( \Gamma_m (X_N ), \Gamma_n( X_N ) ) 
&\stackrel{Thm.\, \ref{mnthm7.5}}{=} 
\mathop{\sum_{1 \leq u \leq m}}_{1 \leq v \leq n}
q_{m,u}' q_{n,v}'
\mathop{\sum_{0 < x \leq u}}_{0 < y \leq v}
\Big| S_{NC} \Big( \myatop{u }{ x}, \myatop{v }{ y} \Big) \Big|_c \\
&\stackrel{Thm. \ref{circularrecursion}}{=}
\mathop{\sum_{0 < x \leq u}}_{0 < y \leq v}
\mathop{\sum_{1 \leq u \leq m}}_{1 \leq v \leq n}
q_{m,u}' q_{n,v}' q_{u,x} q_{v,y}
| S_{NC}( \underline x, \underline y) |_c \\
&=
\mathop{\sum_{0 < x \leq u}}_{0 < y \leq v}
\delta_{m,x} \delta_{n,y}
| S_{NC}( \underline x, \underline y) |_c \\
&=
| S_{NC}( \underline m, \underline n) |_c
\end{align*}
\end{proof}

\noindent
{\bf Proof of Theorem \ref{mt}.} Given that all the
cumulants of degree higher than three for words in $X_{N,1},
\dots , X_{N,k}$ asymptotically vanish, in order to prove
asymptotic independence we only need calculate the first and
second cumulants.

From Proposition \ref{bigprop} we see that $S_{\vec m, \vec
\imath}$ and $S_{\vec n, \vec \jmath}$ are asymptotically
independent unless $(\vec m, \vec \imath)$ and $(\vec n,
\vec \jmath)$ are cyclically equivalent; and that the
complex variance of $S_{\vec m, \vec \imath}$ is $|(\vec m,
\vec \imath)|$.  From Proposition \ref{gsindep} we see that
$\{ \Tr( \Gamma_n( X_{N,i} ) )\}_{n, i}$ are asymptotically
independent from $\{ S_{\vec m, \vec \imath} \}_{\vec m,
\vec \imath}$.  By Proposition \ref{ggindep} we see that
$\Tr( \Gamma_m(X_{N, i} ))$ is asymptotically independent
from $\Tr( \Gamma_n(X_{N, j} ))$ unless $(m,i) = (n, j)$ and
that the variance of $\Tr( \Gamma_m(X_{N,i}) )$ converges to
$m c^m$. What remains is to calculate the asymptotic moments
of $\Tr( \Gamma_n(X_N) )$, $\Tr( \Pi_n(X_N) )$, and $S_{\vec
m, \vec \imath}$.

We start with $\lim_N \Tr( \Pi_n(X_N) )$. Since 
\[
0 = \int_a^b \Pi_n(t) \, d\mu_c(t) = \sum_{k=0}^n p_{n,k}' \sum_{\pi \in
NC(k)} c^{\#(\pi)}
\]
we have
\begin{align*}
\lim_N  \E( \Tr( \Pi_n( X_N ))) &=
\sum_{k=0}^n p_{n,k}' \E( \Tr( X_N^k )) \\
&=
\lim_N \sum_{k=0}^n p_{n,k}' \E( \Tr( X_N^k -   \sum_{\pi \in NC(k)}
c^{\#(\pi)}I_N)) \\
&=
\sum_{k=0}^n p_{n,k}'\sum_{\pi \in NC(k)} c' \#(\pi) c^{\#(\pi) -1}
\end{align*}
Where the last equality is from \cite{mn}, Cor. 9.4.
Hence we must show that
$$
\sum_{k=0}^n p_{n,k}'\sum_{\pi \in NC(k)} \#(\pi) c^{\#(\pi) -1} 
=
\begin{cases}
0 & n \mbox{ even } \\
c^m & n = 2m + 1 
\end{cases}
$$
In terms of our matrix $\Pi$ our claim is that
$$
\Pi 
\begin{pmatrix} 0 \\ 1 \\ 1 + 2 c \\ \vdots \\ \sum_\pi \#(\pi)
c^{\#(\pi) - 1} \\ \vdots \end{pmatrix}
=
\begin{pmatrix} 0 \\ 1 \\ 0 \\ c \\ 0 \\ c^2 \\ \vdots\end{pmatrix}
$$
Or by taking inverses we must show that
$$
\sum_{k=0}^{[\frac{n-1}{2}]} p_{n, 2k+1} c^k =
\sum_{\pi \in NC(n)} \#(\pi) c^{ \#(\pi) -1}
$$
However this is precisely Theorem \ref{lineardecomp}. Thus we have that
$$
\lim_N \E( \Tr( \Pi_n( X_N ) - c \Pi_{n-2}(X_N))) = 0
$$
and thus by equation (\ref{firstsecond})$$
\lim_N \E( \Tr( \Gamma_n( X_N ) + \Gamma_{n-1}(X_N))) = 0
$$
Thus to prove that $\lim_N \E( \Tr( \Gamma_n( X_N )) = (-1)^n c'$ we
check the first few values of $n$ by direct calculation ($n = 0, 1, 2$)
and then obtain the rest by induction using the equation above ($n >
2$). 

Finally we check that $S_{\vec m, \vec \imath}$ is
asymptotically centered. Recall that $p_{k,0} = \sum_{\pi
\in NC(k)} c^{\#(\pi)}$ is the limit of $\E(
\Tr(X_N^k)))$. It is enough to show that $$ \lim_N \E( \Tr(
(X_{N,i_1}^{m_1} - p_{m_1,0} I_N) \cdots (X_{N,i_k}^{m_k} -
p_{m_k,0} I_N) ) ) = 0
$$
and then take linear combinations.

Recall that $\E( \Tr( X_N^m ) ) = N \sum_{\pi \in NC(m)}
\Big(\frac{M}{N}\Big)^{\#(\pi)} + O(N^{-1})$, so it suffices
to keep track only of the terms of order $N$. Thus
\begin{align*}
\lim_N & \E( \Tr( (X_{N,i_1}^{m1} -p_{m_1,0} I_N)
\cdots  (X_{N,i_k}^{mk} -p_{m_k,0} I_N) ) ) \\
&=
\lim_N \sum_{A \subset [k]} \prod_{j \in A} (-p_{m_j,0})
\E( \Tr\Big( \prod_{j \not \in A} X_{N, i_j}^{m_j} \Big) ) \\
&= 
\lim_N N \sum_{A \subset [k]}  \prod_{j \in A} (-p_{m_j,0})
\sum_{\pi \in NC(A^c)} \Big(\frac{M}{N}\Big)^{\#(\pi)}
\end{align*}
were $NC(A)$ is the set of non-crossing partitions of $\cup_{j\in A} [m_{j-1} +
1, m_j]$ that respect the colouring $\vec i$, i.e. $r$ and $s$ are in the
same block of a partition only if $i_r = i_s$. 

Now let us write $NC(A)$ as the disjoint union of the two
subsets $NC(A)_d$ and $NC(A)_c$; where $NC(A)_d$ (\textit{d}
for disconnected) is the set of partitions on $NC(A)$ for
which each block is a subset of an interval $[m_{j-1}+1,
m_j]$ for some $j$, and $NC(A)_c$ is the set of partitions
in $NC(A)$ that connect at least one pair of intervals.
Note that $NC(A)_d = \prod_{j \in A} NC( [m_{j-1}+1,
m_j])$. Let us introduce an additional piece of
notation. For a subset $C \subset [k]$ let $NC(C)_{cc}$ be
the set of non-crossing partitions $\tau$ of $\cup_{j \in C}
[m_{j-1}+1, m_j]$ such that every interval $I_r =
[m_{r-1}+1, m_r]$ is connected to some other interval $I_s =
[m_{s-1}+1, m_s]$ by some block of $\tau$. Thus
$$
NC(A^c)_c = \mathop{\bigcup_{B \cup C = A^c}}_{B \cap C = \emptyset}
NC(B)_d \times NC(C)_{cc}
$$

Note that for every $C \subset [k]$ for which $NC(C)_{cc}$
is not empty we must have $|C^c| \geq 2$ as every
non-crossing partition of $k$ cyclically alternating colours
must have at least two singletons.

Hence, suppressing an $O(N^{-1})$ term, we have 
\begin{align*}
\E( &\Tr( (X_{N,i_1}^{m_1} - p_{m_1,0} I_N) \cdots 
         (X_{N,i_k}^{m_k} - p_{m_k,0} I_N) ) ) \\
&= N \sum_{A \subset [k]}  \prod_{j\in A} (-p_{m_j,0})
\mathop{\sum_{B \cup C = A^c}}_{B \cap C = \emptyset}
\sum_{\pi \in NC(B)_d} \Big(\frac{M}{N}\Big)^{\#( \pi )}
\sum_{\sigma \in NC(C)_{cc}} \Big(\frac{M}{N}\Big)^{\#( \sigma )} \\
&=
N \sum_{C \subset [k]}
\sum_{\sigma \in NC(C)_{cc}} \Big(\frac{M}{N}\Big)^{\#( \sigma )}\kern-10pt
\mathop{\sum_{A \cup B = C^c}}_{A \cap B = \emptyset}
\prod_{j\in A} (-p_{m_j,0})
\prod_{l \in B} \sum_{\sigma \in NC(I_l) }\kern-10pt
\Big(\frac{M}{N}\Big)^{\#( \sigma )} \\
&=
N \sum_{C \subset [k]}
\sum_{\sigma \in NC(C)_{cc}} \Big(\frac{M}{N}\Big)^{\#( \sigma )}
\prod_{j \in C^c} \Big(
\sum_{\pi \in NC(I_j)} \Big(\frac{M}{N}\Big)^{\#( \pi )}
- p_{m_j,0} \Big) 
\end{align*}

Since every $C$ for which $NC(C)_{cc}$ is not empty we have
$|C^c| \geq 2$ there are always at least two factors in the
product in the last expression above. The first factor of
the form
$$
\sum_{\pi \in NC(I)} \Bigg( \frac{M}{N}\Bigg)^{\#(\pi)} - c^{\#(\pi)}
$$ absorbs the $N$ the second converges to 0; thus each term
in the sum above converges to 0 and thus the $S_{\vec m,
\vec \imath}$'s are centered. \qed

\section{Wick Products}\label{wick}

In this section we show how to use non-crossing linear
half-permu\-ta\-tions to obtain relations for the Wick
products of compound Poisson elements. These non-commutative
polynomials were called the free Kailath-Segall polynomials
by M.~Anshelevich \cite{anshelevich}, \S 3.7, which are
in turn a special case of the $q$-Kailath-Segall polynomials
\cite[\S 4.4]{anshelevich}. In Theorem \ref{wickdecomposition}
we show how to expand a monomial as a sum of Wick polynomials.
This is a special case of \cite{anshelevich}, Th. 4.11 (a).
We define the convolution of a
pair of non-crossing linear half-permutations and show in
Theorem \ref{wicktheorem} that this corresponds to the
product of two Wick products. This is  a special case of
\cite{anshelevich}, Thm. 4.11 (c).

We shall use the notation and definitions of \cite{ms},
\S 4.2. Let $\D$ be a unital $*$-algebra equipped with a
tracial state $\psi$ and represent $\D$, via the
GNS-representation, on $\HH:=
\overline{\D}^{\la\cdot,\cdot\ra}$, where the inner product on
$\HH$ is given by $$\la d_1,d_2\ra:=\psi(d_2^*d_1).$$ Let
$\FF(\HH)$ be the full Fock space over $\HH$ with vacuum state
$\Omega$ and for $d \in \D$ we have the following operators on
$\FF(\HH)$: $l(d)$ is the left creation operator, $l^\ast(d)$
is the left annihilation operator, $\Lambda(d)(d_1 \ot \cdots
\ot d_n) = (d\, d_1) \ot d_2 \ot \cdots d_n$ is the
preservation operator (with $\Lambda(d) \Omega = 0$), and
$p(d) = l(d) + l^\ast(d^\ast) + \Lambda(d) + \psi(d)1$.
$\Omega$ is a cyclic and separating vector for the algebra
generated by $\{ p(d) \mid d \in \D\}$. Thus for each $d_1 \ot
\cdots \ot d_n \in \FF(\HH)$ there is a unique polynomial
$W(d_1 \ot \cdots \ot d_n)$ in the non-commuting variables $\{
p(d) \mid d \in \D\}$ such that $W(d_1 \ot \cdots \ot
d_n)\Omega = d_1 \ot \cdots \ot d_n$. These are the {\it Wick
products}.

One checks that $W(d) = p(d) - \psi(d)$ and from the
definition of $p(d)$ on has immediately that
\begin{multline}\label{wickdef}
W(d \ot d_1 \ot \cdots \ot d_n) = p(d) W(d_1 \ot \cdots \ot
d_n) \\ - \psi(d\, d_1) W(d_2 \ot \cdots \ot d_n) - W((d\,
d_1) \ot \cdots \ot d_n) \\ - \psi(d) W(d_1 \ot d_2 \ot
\cdots \ot d_n)
\end{multline}

Our first goal in this section is to prove the following
\begin{equation}\label{reverse}
W(d_1 \ot \cdots \ot d_n)^\ast = W(d_n^\ast \ot \cdots \ot
d_1^\ast)
\end{equation}
In the course of proving this we shall extend the definition of  Wick products
to $W_\pi(d_1 \ot \cdots \ot d_n)$ where $\pi$ is a non-crossing linear
half-permutation on $[n]$. When $\pi$ is the half-permutation in which each
block is an open singleton $W_\pi(d_1 \ot \cdots \ot d_n) = W(d_1 \ot \cdots
\ot d_n)$. 

\begin{definition}
Let $(\pi, B_1, \dots , B_k)$ be a non-crossing linear
half-permu\-tation on $[k]$ with open blocks $B_1, \dots ,
B_k$ and closed blocks $C_1, \dots , C_l$.  For a cycle $C =
(r_1, \dots , r_u)$ let $\psi_C(d_1, \dots , d_n) =
\psi(d_{r_1} \cdots d_{r_u})$. Write each cycle $B_j =
(r_{j,1}, \dots r_{j, v(j)} )$. Then we define
\begin{multline*}
W_\pi( d_1 \ot \cdots \ot d_n) = 
\psi_{C_1}(d_1, \dots , d_n) \cdots \psi_{C_l}(d_1, \dots , d_n) \\
\times W\big( (d_{r_{1,1}} \cdots d_{r_{1, v(1)}}) \ot \cdots \ot
(d_{r_{k,1}} \cdots d_{r_{k,v(k)}} ) \big)
\end{multline*}
Finally suppose that $S$ is a finite set of non-crossing
linear half-permutations. Let $$W_S(d_1 \ot \cdots \ot d_n)
= \sum_{\sigma \in S} W_\sigma(d_1 \ot \cdots \ot d_n)$$
\end{definition}

\begin{example}
Let $\pi = (1, 2) (3) (4) (5,6)$ with $(1, 2)$ and $(4)$
open blocks and $(3)$ and $(5, 6)$ closed. Then
$$W_\pi(d_1 \ot \cdots \ot d_6) = \psi(d_3) \psi(d_5\,
d_6) W\big( (d_1 d_2) \ot d_4 \big)$$
\end{example}

\begin{theorem}\label{wickdecomposition}
Let $d_1, d_2, \dots d_n \in D$ and $\pi \in NCL(n)$. Then
\begin{enumerate}
\item for $d_0 \in \D$
\begin{multline*}
p(d_0) W_\pi( d_1 \ot \cdots \ot d_n) =
W_{\dot\pi}(d_0 \ot  d_1 \ot \cdots \ot d_n) 
+ W_{\ddot\pi}(d_0 \ot  d_1 \ot \cdots \ot d_n) \\ + 
W_{\dot{\ddot\pi}}(d_0 \ot  d_1 \ot \cdots \ot d_n)   +
W_{\ddot{\ddot\pi}} (d_0 \ot d_1 \ot \cdots \ot d_n)
\end{multline*}
where the half-permutations $\dot\pi$, $\ddot\pi$,
$\dot{\ddot\pi}$, $\ddot{\ddot\pi}$ of $\{0 \} \cup [n]$ are
obtained as follows:  \begin{itemize}
\item $\dot\pi$: by adding an open singleton at 0;
\item $\ddot\pi$: by joining 0 to the first open block and then making this a
closed block;
\item $\dot{\ddot\pi}$:  by joining 0 to the first open block and
leaving it open; 
\item $\ddot{\ddot\pi}$: by making 0 a closed
singleton. \end{itemize}
\item $p(d_1) \cdots p(d_n) =
\ds\sum_{\pi \in \ncl(n)} W_\pi( d_1 \ot \cdots \ot d_n)$;
\end{enumerate}

\end{theorem}

\begin{proof}
({\it a}) follows from equation (\ref{wickdef}) when we rewrite it as
\begin{multline*}
p(d_0) W(d_1 \ot \cdots \ot d_n) 
=
 W(d_0 \ot d_1 \ot \cdots \ot d_n) \\
\mbox{} + \psi(d_0\, d_1) W(d_2 \ot \cdots \ot d_n)
\mbox{} + W((d_0d_1) \ot \cdots \ot d_n) \\  
\mbox{} + \psi(d_0) W(d_1 \ot d_2 \ot \cdots \ot d_n)
\end{multline*}

({\it b}) follows by induction and the fact that every
$\sigma \in \ncl(\{0 \} \cup [n])$ is $\dot\pi$, $\ddot\pi$,
$\dot{\ddot\pi}$, or $\ddot{\ddot\pi}$ for exactly one $\pi
\in \ncl([n])$. \end{proof}

\begin{corollary}\label{wickadjoints}
For $d_1 \dots , d_n \in D$, $$W(d_1 \ot \cdots \ot
d_n)^\ast = W(d_n^\ast \ot \cdots \ot d_1^\ast)$$
\end{corollary}

\begin{proof}
Take adjoints of equation ({\it b}) in Theorem
\ref{wickdecomposition} and use induction and the fact that
$\ncl(n)$ is invariant under the reflection $i \mapsto n+1
-i$.  \end{proof}

\begin{notation}
Let $\pi \in \ncl(m)$ and $\sigma \in \ncl(n)$ be
non-crossing linear half-permutations. Let the open blocks
of $\pi$ be $B_1, \dots , B_j$ and the open blocks of
$\sigma$ be $C_1, \dots , C_k$ and $l = \min\{j, k\}$. We
shall regard $\sigma$ as a permutation of $[m+1, m+n]$. We
shall construct $2l + 1$ non-crossing half-permutations $\pi
\vee_0 \sigma, \pi \vee_{1,o} \sigma, \pi \vee_{1,c} \sigma,
\dots, \pi \vee_{l,o} \sigma, \pi \vee_{l,c} \sigma$ in
$\ncl(m+n)$. The set $\{ \pi \vee_0 \sigma , \dots, \pi
\vee_{l,c} \sigma \}$ will be denoted $\pi \ast \sigma$, and
called the \textit{convolution} of $\pi$ and $\sigma$.

Let $\pi \vee_0 \sigma$ be the permutation of $[m + n]$
obtained by letting $\pi$ act on $[m]$ and $\sigma$ act on
$[m + 1, m + n]$.

\begin{center}\leavevmode
\hbox{\vbox{\hsize115pt\raggedright\parindent0pt
\hbox{\includegraphics{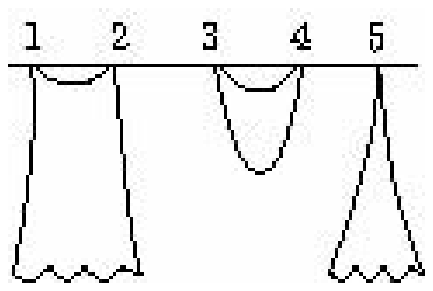}}\leavevmode
\hbox{}\hfill$\pi$\hfill\hbox{}}\qquad
\vbox{\hsize155pt\raggedright\parindent0pt
\hbox{\includegraphics{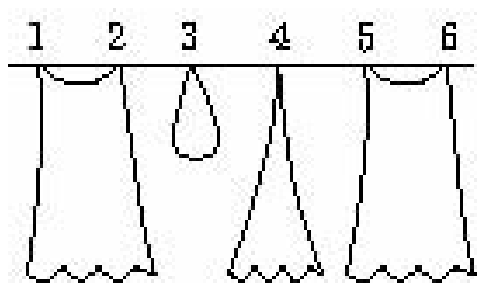}}\leavevmode
\hbox{}\hfill$\sigma$\hfill\hbox{}}}
\end{center}
\medskip

\begin{center}\leavevmode \hbox{\includegraphics{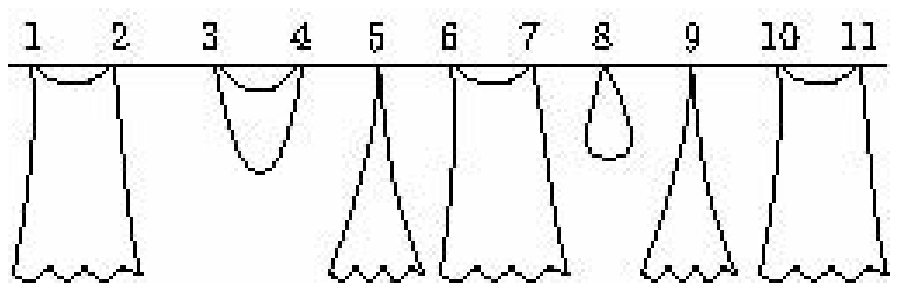}}
\medskip

\leavevmode\vbox{\hsize 250pt\noindent\raggedright\small
\textbf{Figure \figno.} $\pi \vee_0 \sigma$ is the concatenation
of the two non crossing linear half-permutations.}

\end{center}\medskip

In the figure above we have $\pi = \{ (1, 2), (3,4), (5)\}$
with (1, 2) and (5) open blocks; and $\sigma = \{(1, 2),
(3), (4),(5,6)\}$ with (1, 2) , (4), and (5, 6) open
blocks. We transfer $\sigma$ to the interval [6, 11] and
obtain $\pi \vee_0 \sigma = \{ (1, 2), (3,4),
(5),(6,7),(8),(9), (10, 11)\}$ with (1,2), (5), (6,7), (9),
and (10,11) open.

The open blocks of $\pi$ and $\sigma$ are ordered by their
smallest element, which is the same as ordering by their
largest element, as open blocks cannot be nested. The basic
construction is to take the rightmost open block of $\pi$,
$B_j = (p_1, \dots , p_r)$ and the leftmost open block of
$\sigma$, $C_1 =(q_1, \dots , q_s)$, and replace them with
$(p_1, \dots , p_r, q_1, \dots , q_s)$ first as an open
block then as a closed block. These permutations are called
$\pi \vee_{1,o} \sigma$ and $\pi \vee_{1,c} \sigma$
respectively. We obtain $\pi \vee_{2,o} \sigma$ and $\pi
\vee_{2,c} \sigma$ by applying this procedure (after $B_j$
and $C_1$ have been joined as a closed block) to the open
blocks $B_{j-1}$ and $C_{2}$ and continuing until all the
open blocks of one of the permutations are exhausted.  This
will produce $2l +1$ permutations.  \end{notation}

\begin{eqnarray*}
\raise-4pt\hbox{\includegraphics{fig20.eps}\ } \ast 
\raise-4pt\hbox{\ \includegraphics{fig20.eps}} & =&
\raise-4pt\hbox{\includegraphics{fig21.eps}\ } +
\raise-4pt\hbox{\ \includegraphics{fig22.eps}\ } +
\raise-4pt\hbox{\ \includegraphics{fig23.eps}} \\ [10pt]
&& \mbox{} +
\raise-4pt\hbox{\ \includegraphics{fig24.eps}\ } +
\raise-4pt\hbox{\ \includegraphics{fig25.eps}}
\end{eqnarray*}

\begin{lemma}\label{inductivewick}
$$\displaylines{
W(d_1 \ot \cdots \ot d_m)(e_1 \ot \cdots \ot e_n) 
\hfill \cr \hfill \mbox{} - \psi(d_m e_1)
W(d_1 \ot \cdots \ot d_{m-1})(e_2 \ot \cdots \ot e_n)  \hfill\cr
\mbox{} =
d_1 \ot \cdots \ot d_m \ot e_1 \ot \cdots \ot e_n \hfill\cr
\hfill \mbox{} +
d_1 \ot \cdots \ot d_m e_1 \ot \cdots \ot e_n \hfill 
\cr}$$
\end{lemma}

\begin{proof}
We prove the lemma by induction on $n$. By writing 
$$\displaylines{
e_1 \ot \cdots \ot e_n =
p(e_1) e_2 \ot \cdots \ot e_n - e_1 e_2 \ot \cdots \ot e_n \cr
\mbox{} -
\psi(e_1 e_2) e_3 \ot \cdots \ot e_n - \psi(e_1) e_2 \ot \cdots
\ot e_n \cr}$$
taking the adjoint of equation \ref{wickdef} and combining
this with Corollary \ref{wickadjoints} we get that
$$\displaylines{
W(d_1 \ot \cdots \ot d_m)(e_1 \ot \cdots \ot e_n) 
\hfill \cr \hfill \mbox{} - \psi(d_m e_1)
W(d_1 \ot \cdots \ot d_{m-1})(e_2 \ot \cdots \ot e_n)  \hfill\cr
\mbox{} =
W(d_1 \ot \cdots \ot d_m \ot e_1)(e_2 \ot \cdots \ot e_n) 
\hfill \cr \hfill \mbox{} - \psi(e_1 e_2)
W(d_1 \ot \cdots \ot d_{m})(e_3 \ot \cdots \ot e_n)  \hfill\cr
\hfill\mbox{} +
W(d_1 \ot \cdots \ot d_m e_1)(e_2 \ot \cdots \ot e_n) \hfill\cr
\hfill \mbox{} -
W(d_1 \ot \cdots \ot d_m)(e_1 e_2 \ot \cdots \ot e_n) 
\cr}$$

We may write
\begin{align*}
W(d_1 &  \ot \cdots \ot d_m e_1)(e_2 \ot \cdots \ot e_n)  -
W(d_1 \ot \cdots \ot d_m)(e_1 e_2 \ot \cdots \ot e_n) \\
&=
\Big\{ W(d_1 \ot \cdots \ot d_m e_1)(e_2 \ot \cdots \ot e_n)  \\
& \qquad - \psi(d_m e_1 e_2) W(d_1 \ot \cdots \ot d_{m-1})(e_3 \ot
\cdots \ot e_n) \Big\} \\
&\quad\mbox{} -
\Big\{ W(d_1 \ot \cdots \ot d_m)(e_1 e_2 \ot \cdots \ot e_n) \\
& \qquad \mbox{} - \psi(d_m e_1
e_2) W(d_1 \ot \cdots \ot d_{m-1})(e_3 \ot \cdots \ot e_n) \Big\}
\end{align*}

By applying our induction hypothesis to the bracketed terms
we get
\begin{align*} 
&W(d_1 \ot \cdots \ot d_m)(e_1 \ot \cdots \ot e_n) \\
&\qquad - \psi(d_m e_1)
W(d_1 \ot \cdots \ot d_{m-1})(e_2 \ot \cdots \ot e_n)  \\
& \mbox{} =
W(d_1 \ot \cdots \ot d_m \ot e_1)(e_2 \ot \cdots \ot e_n) \\
&\quad\mbox{} - \psi(e_1 e_2)
W(d_1 \ot \cdots \ot d_{m})(e_3 \ot \cdots \ot e_n)  \\
& \quad\mbox{} +
d_1  \ot \cdots \ot d_m e_1 \ot \cdots \ot e_n 
 -
d_1 \ot \cdots \ot d_m \ot e_1 e_2 \ot \cdots \ot e_n
\end{align*}
To conclude the proof we apply this formula $n-1$
times to obtain
$$\displaylines{
W(d_1 \ot \cdots \ot d_m)(e_1 \ot \cdots \ot e_n) 
\hfill \cr \hfill \mbox{} - \psi(d_m e_1)
W(d_1 \ot \cdots \ot d_{m-1})(e_2 \ot \cdots \ot e_n)  \hfill\cr
\mbox{} =
W(d_1 \ot \cdots \ot e_n)\Omega + d_1 \ot \cdots \ot d_m e_1 \ot
\cdots \ot e_n
\cr}$$
as required.
\end{proof}

\begin{theorem}\label{wicktheorem}
$$
W_\pi(d_1 \ot \cdots \ot d_m) W_\sigma(e_1 \ot \cdots \ot e_n) =
W_{\pi \ast \sigma}(d_1 \ot \cdots \ot d_m \ot e_1 \ot \cdots
\ot e_n)
$$ 
\end{theorem}

\begin{proof}From the definition we only have to verify this in
the special case that all the blocks of $\pi$ and $\sigma$
are open singletons. From the lemma above we have

$$\displaylines{
W(d_1 \ot \cdots \ot d_m) W(e_1 \ot \cdots \ot e_n) \hfill\cr
\mbox{} =
d_1 \ot \cdots \ot d_m \ot e_1 \ot \cdots \ot e_n 
+
d_1 \ot \cdots \ot d_m e_1 \ot \cdots \ot e_n \cr
\mbox{} +
\psi(d_m e_1) W(d_1 \ot \cdots \ot d_{m-1})
W(e_2 \ot \cdots \ot e_n)
\cr}$$
Thus the theorem follows from repeated applications of lemma
\ref{inductivewick} above.
\end{proof}

\section{The case of pairings}

In this paper we have only dealt with Wishart matrices, however
there is a parallel program one can carry out for Gaussian
matrices. Indeed Theorem \ref{mt} in the Gaussian case was
obtained by Cabanal-Duvillard \cite{thierry} using stochastic
integration. However a combinatorial proof along the lines
presented here is equally possible. The idea is to consider only
pairings throughout all the calculations. One has the notion of
non-crossing circular and linear half-pairings; in this case the
polynomials are the Chebyshev polynomials $\{ C_n \}_n$ and
$\{S_n \}_n$ of the first and second kind respectively
renormalized to the interval $[-2, 2]$. Exactly as in Corollary
\ref{wickadjoints} one obtains the analogous result on Wick
products needed in \cite{ms}, Lemma 5.5. Theorem
\ref{wicktheorem} likewise has a similar formulation and proof.
The pairing version of some of our results (e.g. Theorem
\ref{wickdecomposition}) are special cases of Effros and Popa
\cite{ep}, Thm 3.3.

\section*{Acknowledgments}
The authors wish to thank Alexandru Nica for fruitful
discussions at an early stage of this work.

%\phantom{X}\footnotetext{{last revised \today}}


\begin{thebibliography}{VVI}

\bibitem[{\sc a}]{anshelevich} {\it M. Anshelevich}, Appel
Polynomials and their relatives, 
Inter. Math. Res. Notes, 2004, no. 65, 3469-3531.

\bibitem[{\sc b}]{biane} {\it P. Biane}, Some properties of
crossings and partitions, Discrete Mathematics, {\bf
175}, (1997), 41-53.

\bibitem[{\sc cd}]{thierry} {\it T. Cabanal-Duvillard},
Fluctuations de la Loi Empirique de Grandes Matrices
Al\'eatoires, Ann. I. H. Poincar\'e (B), Probabilit\'es
et Statistiques, {\bf 37} (2001) 373 - 402.

\bibitem[{\sc ep}]{ep} {\it E. G. Effros and M. Popa}, Feynman
diagrams and Wick Products associated with $q$-Fock space,
Proc. Nat. Acad. Sci. (USA), {\bf 100}, 8629 -
8633.

\bibitem[{\sc joh}]{kurt} {\it K. Johansson}, On Fluctuations of
Eigenvalues of Random Hermitian Matrices, Duke
Math. J., {\bf 91} ( 1998), 151 - 204.

\bibitem[{\sc jon}]{jon} {\it D. Jonsson} Some limit theorems
for the eigenvalues of a sample covariance matrix,
J. Mult. Anal., {\bf 12} (1982), 1-38.

\bibitem[{\sc ht$_1$}]{ht1} {\it U. Haagerup and
S. Thorbj{\o}rnsen}, Random Matrices and $K$-theory for Exact
C*-algebras, Documenta Mathematica, {\bf 4} (1999),
341 - 450.

\bibitem[{\sc ht$_2$}]{ht2} {\it Uffe Haagerup and
S. Thorbj{\o}rnsen}, Random Matrices with Complex Gaussian
Entries, Expositiones Math., {\bf 21} (2003), 293 -
337.

\bibitem[{\sc hp}]{hp} {\it F. Hiai and D. Petz}, The
Semicircle Law, Free Random Variables and Entropy,
Providence, R.I., Amer. Math.  Soc., 2000.

\bibitem[{\sc j}]{jones} {\it V. F. R. Jones}, The annular
structure of subfactors, L'Enseignement Mathematique,
to appear, math.OA/0105071

\bibitem[{\sc mn}]{mn} {\it J. Mingo and A. Nica}, Annular
non-crossing permutations and partitions, and second-order
asymptotics for random matrices, 
Inter. Math. Res. Notices, 2004, n$^{\rm o}$ 28, 1413 -
1460.

\bibitem[{\sc ms}]{ms} {\it J. Mingo and R. Speicher}, Second
order freeness and fluctuations of random matrices:
I. Gaussian and Wishart matrices and cyclic Fock spaces,
Preprint, 2004, math/OA/0405191.

\bibitem[{\sc m\'ss}]{mss} {\it J. Mingo, P. \'Sniady, and
R. Speicher} Second Order Freeness and Fluctuations of
Random Matrices: II. Unitary Random Matrices, preprint May
2004, 22 pages. math/OA/0405258


\end{thebibliography}
\end{document}